\documentclass[10pt,a4paper]{amsart}

\usepackage[T1]{fontenc}
\usepackage[utf8]{inputenc}
\usepackage{lmodern} 
\usepackage[english]{babel}
\usepackage{csquotes}
\usepackage[a4paper,margin=1in,marginparwidth=2cm]{geometry}
\usepackage{microtype}
\usepackage{xcolor}
\usepackage{graphicx}
\usepackage{amsmath,amssymb,amsthm,amsfonts,mathrsfs,amsopn}
\usepackage{mathtools}
\mathtoolsset{showonlyrefs} 
\usepackage{esint,slashed}
\usepackage[all]{xy}
\usepackage{tikz-cd}
\usetikzlibrary{arrows.meta}
\definecolor{ladderblue}{RGB}{30,72,125}
\definecolor{laddergray}{RGB}{80,87,98}

\newif\ifdraftnotes
\draftnotesfalse
\newif\ifshowreferencekeys
\showreferencekeysfalse
\ifdraftnotes
  \usepackage[textsize=footnotesize]{todonotes}
\else
  \usepackage[disable]{todonotes}
\fi

\ifshowreferencekeys
  \usepackage[notref,notcite]{showkeys}
\fi

\usepackage[
  backend=biber,
  style=numeric-comp,
  sorting=nyt,
  sortcites=true,
  giveninits=true,
  maxbibnames=99,
  doi=true,
  url=true,
  eprint=true,
  isbn=false
]{biblatex}
\DeclareFieldFormat[article,inproceedings,unpublished,misc,online]{title}{#1}
\renewbibmacro{in:}{\ifentrytype{article}{}{\printtext{\bibstring{in}\intitlepunct}}}

\usepackage{xpatch}
\makeatletter
\xpatchbibmacro{cite:dump}
  {\ifnumgreater{\value{cbx@tempcnta}}{1}}
  {\ifnumgreater{\value{cbx@tempcnta}}{0}}
  {}
  {\PackageError{math-paper-template}{Citation compression patch failed}
    {Check the numeric-comp citation style in your TeX installation.}}
\makeatother

\AtEveryBibitem{
  \iffieldundef{doi}
    {\iffieldundef{eprint}
       {}
       {\clearfield{url}\clearfield{urlyear}\clearfield{urlmonth}\clearfield{urlday}}}
    {\clearfield{url}\clearfield{eprint}
     \clearfield{urlyear}\clearfield{urlmonth}\clearfield{urlday}}
}

\usepackage{xurl}
\usepackage{hyperref}
\definecolor{linkblue}{RGB}{25,55,105}
\hypersetup{
  colorlinks=true,
  linkcolor=linkblue,
  citecolor=linkblue,
  urlcolor=linkblue,
  pdfborder={0 0 0},
  pdftitle={Dirac eigenvalue comparisons},
  pdfauthor={}
}
\usepackage{bookmark}
\DeclareFieldFormat{doi}{%
  \mkbibacro{DOI}\addcolon\space
  \href{https://doi.org/#1}{\nolinkurl{#1}}}
\DeclareFieldFormat{eprint:arxiv}{%
  \newline\href{https://arxiv.org/abs/#1}{\nolinkurl{https://arxiv.org/abs/#1}}}
\DeclareFieldAlias{eprint:arXiv}{eprint:arxiv}
\DeclareFieldFormat{url}{\url{#1}}

\newcommand{\R}{\mathbb{R}}
\newcommand{\C}{\mathbb{C}}
\newcommand{\bZ}{\mathbb{Z}}

\newcommand{\cO}{\mathcal{O}}

\newcommand{\dd}{\mathop{}\!\mathrm{d}}
\newcommand{\eps}{\varepsilon}
\newcommand{\p}{\partial}
\newcommand{\sph}{\mathbb{S}}
\newcommand{\tor}{\mathbb{T}}

\newcommand{\D}{\mathbin{\slashed{D}}}

\newcommand{\Penr}{\slashed{P}}
\newcommand{\snabla}{\slashed{\nabla}} 

\DeclareMathOperator{\dv}{\dd{vol}}

\DeclareMathOperator{\diverg}{div}

\DeclareMathOperator{\grad}{grad}

\DeclareMathOperator{\id}{Id}
\DeclareMathOperator{\Ker}{Ker}

\DeclareMathOperator{\Ric}{Ric}

\DeclareMathOperator{\Spec}{Spect}

\DeclareMathOperator{\vol}{vol}

\theoremstyle{plain}
\newtheorem{thm}{Theorem}[section]
\newtheorem{lemma}[thm]{Lemma}
\newtheorem{prop}[thm]{Proposition}
\newtheorem{cor}[thm]{Corollary}

\newtheorem{theorem}{Theorem}

\theoremstyle{definition}

\theoremstyle{remark}
\newtheorem{rmk}[thm]{Remark}

\title[Dirac eigenvalue comparison]{Dirac eigenvalue comparison}

\thanks{P.S. and R.W. are partially supported by the Fundamental Research Funds for the Central Universities and Beijing Institute of Technology Research Fund Program for Young Scholars.
R.W. is partially supported by Beijing Natural Science Foundation No. 1262018.}

\author[Z. Qiu]{Zixuan Qiu}
\address{Zixuan Qiu, School of mathematics and statistics, Beijing Institute of Technology, Zhongguancun South Street No. 5, 100081 Beijing, P.R.China.}
\email{ZixuanQiu@bit.edu.cn}
\author[P. Shi]{Pengshuai Shi}
\address{Pengshuai Shi, School of mathematics and statistics, Beijing Institute of Technology, Zhongguancun South Street No. 5, 100081 Beijing, P.R.China.}
\email{shipengshuai@bit.edu.cn}
\author[R. Wu]{Ruijun Wu}
\address{Ruijun Wu, School of mathematics and statistics, Beijing Institute of Technology, Zhongguancun South Street No. 5, 100081 Beijing, P.R.China.}
\email{ruijun.wu@bit.edu.cn}
\date{}

\begin{document}
\begin{abstract}
We study how curvature bounds control the Dirac eigenvalues of closed spin manifolds, comparing the positive and negative eigenvalues separately.
Within a fixed conformal class and spin structure, pointwise ordering of scalar curvature relative to a metric of negative scalar curvature implies comparison at every nonzero index.
On the two-sphere, a holomorphic ladder on the half-canonical bundle proves the corresponding comparison in both directions under a Gaussian curvature bound.
On surfaces of higher genus, the same bundle identities propagate comparison from the scalar and spinor levels to higher canonical powers.
In the round conformal class on the higher dimensional sphere, a prescribed scalar curvature bound yields an upper bound in absolute value for the first~$2^{[ n/2 ]}$ Dirac eigenvalues.
When the conformal class varies, we can show that the corresponding comparisons on hyperbolic surfaces and odd-dimensional spheres no longer hold;
while on even-dimensional spheres, the first gap and all sufficiently high-index eigenvalues of either sign satisfy a one-direction comparison.
\end{abstract}
\maketitle

\section{Introduction}

Estimates for Dirac eigenvalues on a compact spin Riemannian manifold are among the central topics in spin differential geometry.
Positive scalar curvature excludes harmonic spinors by Lichnerowicz's formula~\cite{Lichnerowicz1963} and yields a positive lower bound for the absolute values of Dirac eigenvalues by Friedrich's estimate~\cite{Friedrich1980erste}.
Hijazi sharpened Friedrich's estimate using the conformal Laplacian~\cite{Hijazi1986}.
Under appropriate positivity assumptions, Hijazi and Raulot obtained further lower bounds from the Paneitz operator and Branson's~$Q$-curvature~\cite{HijaziRaulot2007Branson}.
Volume-normalized Dirac eigenvalues also lead to conformal spectral invariants, beginning with the first positive eigenvalue studied by Ammann in~\cite{Ammann2003SpinConformal}.

When scalar curvature is negative, its contribution to the Schr\"odinger--Lichnerowicz formula no longer provides a positive lower bound.
A different question is nevertheless natural: once a conformal class and a spin structure have been fixed, does curvature ordering imply ordering of the Dirac eigenvalues?
The conformal identification of metric spin structures allows this question to be formulated on one fixed spinor bundle~\cite{Hitchin1974harmonic,Ammann2001spectral}.
The dimension of the Dirac kernel is constant within a conformal class~\cite{Hitchin1974harmonic}, so we can omit zero and order the eigenvalues on each side of it.
If the conformal class of a higher-genus surface varies, the kernel dimension may jump even when the spin structure is fixed; see Section~\ref{sec:varying-conformal-class}.
Thus we cannot expect a general all-index eigenvalue comparisons from curvature comparison in the full space of Riemannian metrics.

We distinguish the signed Dirac spectrum from the spectrum of its square.
For~$k\in\mathbb Z_*=\mathbb Z\setminus\{0\}$, let~$\lambda_k(g)$ denote the nonzero eigenvalues of~$\D_g$, counted with complex multiplicity and ordered as follows:
\begin{align}
 -\infty\leftarrow\cdots \le\lambda_{-2}(g)\le\lambda_{-1}(g)<0<\lambda_1(g)\le\lambda_2(g)\le\cdots\rightarrow +\infty.
\end{align}
Ammann's conformal estimate~\cite[Proposition 5.2]{Ammann2001spectral} compares the ordered eigenvalues of~$\D^2$.
In the absence of harmonic spinors its first value is~$\min\{(\lambda_1)^2,(\lambda_{-1})^2\}$; in an asymmetric spectrum this does not separate~$\lambda_1$ from~$\lambda_{-1}$.
The variational characterization and the rigidity for scalar weights in~\cite{QiuWu2025} allow us to compare the eigenvalues of either sign, even in the presence of harmonic spinors.
Combining this comparison with the maximum principle for the conformal scalar-curvature equation, we obtain the following comparison with rigidity.

\begin{theorem}\label{thm:negative-intro}
Let~$(M^n,h)$ be a closed spin manifold,~$n\ge2$, with negative scalar curvature~$R_h<0$.
Fix the topological spin structure and let~$g\in[h]$.
Then, for every~$k\ne0$,
\begin{align}
 R_g\ge R_h&\quad\Longrightarrow\quad|\lambda_k(g)|\le|\lambda_k(h)|,\\
 R_g\le R_h&\quad\Longrightarrow\quad|\lambda_k(g)|\ge|\lambda_k(h)|.
\end{align}
In either case, equality for some~$\lambda_k$ holds if and only if~$g=h$.
\end{theorem}
Theorem~\ref{thm:negative-intro} is proved as Theorem~\ref{thm:negative-reference} below.
The scalar curvature~$R_g$ is allowed to change sign in the first implication.
If~$R_h=-c<0$ is constant, we may rescale the reference metric and conclude that
\begin{align}
 -cb^2\le R_g\le-ca^2,\quad 0<a\le b
 \quad\Longrightarrow\quad
 a|\lambda_k(h)|\le|\lambda_k(g)|\le b|\lambda_k(h)|\qquad(k\ne0),
\end{align}
with homothety rigidity.
In particular, these inequalities apply to a fixed hyperbolic conformal class in every dimension.

The maximum principle used in Theorem~\ref{thm:negative-intro} has the wrong sign for a positively curved reference.
A M\"obius transformation of the Riemann sphere that is not a round isometry preserves the Gaussian curvature of the round metric under pullback, although its conformal factor takes values both above and below one.
Thus curvature equality does not imply pointwise ordering of conformal metrics.

Lin--Wang--Xu introduced a holomorphic ladder to compare every positive eigenvalue of the scalar Laplacian on a two-sphere with~$K_g\ge1$ to its round counterpart, with rigidity~\cite{LWX2026}.
Ding--Ge--Li used this ladder to treat both curvature directions, including~$K_g\le1$~\cite{DingGeLi2026}.
Starting from functions, the scalar ladder passes through the anticanonical powers~$\mathcal K^{-m}$, where~$\mathcal K$ denotes the canonical line bundle of~$\mathbb P^1$.
For the Dirac operator we start instead with the positive half-spinor bundle~$S$, a square root of~$\mathcal K$.
Tensoring with~$\mathcal K^{-1}$ determines the successive bundle classes in the Picard group of~$\mathbb P^1$, while the Dolbeault operators act between their spaces of smooth sections after metric identification.
This change in the initial bundle shifts both the eigenvalues and the multiplicities in the ladder.

\begin{theorem}\label{thm:sphere-intro}
Let~$g_0$ be the unit round metric on~$\sph^2$.
For any smooth metric~$g$ with Gaussian curvature~$K_g$ and every~$k\ne0$,
\begin{align}
 K_g\ge1&\quad\Longrightarrow\quad|\lambda_k(g)|\ge|\lambda_k(g_0)|,\\
 K_g\le1&\quad\Longrightarrow\quad|\lambda_k(g)|\le|\lambda_k(g_0)|.
\end{align}
In either case, equality for some~$\lambda_k$ holds if and only if~$(\sph^2,g)$ is isometric to the unit round sphere.
\end{theorem}
Theorem~\ref{thm:sphere-intro} is proved as Theorem~\ref{thm:sphere} below.
Note that no lower curvature assumption is required for the second implication.

Chirality pairs the positive and negative eigenvalues in dimension two.
In Section~\ref{subsec:sphere} we compute the local adjoints and follow the spectral index shifts.
On the unit two-sphere the successive comparisons become equalities, and we recover the eigenvalues~$\pm j$, each of complex multiplicity~$2j$, first computed by S.~Sulanke~\cite{Sulanke1981Berechnung}; see also~\cite[Theorem~2.1.3]{Ginoux2009DiracSpectrum}.

On a surface of genus greater than one with a fixed conformal and spin structure, the maximum principle compares the scalar and spinor levels relative to a hyperbolic reference metric.
The bundle identities then propagate the comparison to all positive integer and half-integer powers of~$\mathcal K$.
Although we cannot compute the Dirac eigenvalues of the hyperbolic reference metric~$g_0$ explicitly, the method is robust enough to give a comparison of the Dirac eigenvalues provided that we have a comparison for the Gaussian curvatures.
Section~\ref{subsec:hyperbolic-ladder} contains this Dolbeault proof and the rigidity statement for the higher twists.

It is natural to ask whether an upper sectional-curvature bound suffices to give eigenvalue comparison on higher-dimensional spheres when the conformal class is allowed to vary.
This turns out to be wrong: already on~$\sph^3$, the first Dirac gap can be arbitrarily large, by shrinking the circle fibers in the Hopf fibration~$\sph^1\hookrightarrow\sph^3\to \sph^2$.
More generally, we get the following no-comparison result. 

\begin{theorem}\label{thm:odd-intro}
For every odd~$n\ge3$, there is a family of smooth metrics~$g_r$ on~$\sph^n$, indexed by~$0<r\le r_0$, such that
\begin{align}
 0<\sec_{g_r}\le1,\qquad
 \min_{\lambda\in\Spec(\D_{g_r})}|\lambda|\longrightarrow\infty\qquad\text{as }r\longrightarrow0.
\end{align}
\end{theorem}
Theorem~\ref{thm:odd-intro} is proved as Theorem~\ref{thm:odd-spheres} below.
The underlying mechanism is the nonprojectable circle collapse developed by Ammann and Ammann--B\"ar~\cite{Ammann1998,AmmannBaer1998}.
For the ordinary Hopf action on~$\sph^{2m+1}$, one full turn acts on spinors by~$(-1)^{m+1}$.
In dimension~$4\ell+1$ the lift therefore excludes zero fiber frequency, whereas in dimension~$4\ell+3$ it allows this frequency and some eigenvalues remain bounded during collapse.
We choose a weighted Sasakian action with odd total weight: one full turn then acts by~$-1$ in every odd dimension, and the same collapse estimate applies while the sectional curvatures remain positive and uniformly bounded above.

In even dimensions the same curvature assumptions lead to an upper bound for the first gap by Klingenberg's injectivity-radius theorem and B\"ar's hyperspherical-radius estimate~\cite{Klingenberg1959,Baer2026Hyperspherical}.
We reprove this consequence and its rigidity for completeness.
There is also a comparison at the other end of the spectrum: if~$0<\sec_g\le1$ on an even-dimensional sphere and~$g$ is not isometric to the unit round metric~$g_0$, then
\begin{align}\label{eq:round-tail-intro}
 |\lambda_k(g)|<|\lambda_k(g_0)|\qquad\text{for all sufficiently large }|k|.
\end{align}
Here Klingenberg's injectivity-radius bound and Berger's volume inequality imply strict volume comparison for a nonround metric; see Corollary~\ref{cor:even-sphere-high-indices}.
The threshold may depend on~$g$.
For the opposite curvature bound~$\sec_g\ge1$, Friedrich's theorem already implies the lower bound for the first gap under the weaker hypothesis~$R_g\ge n(n-1)$~\cite{Friedrich1980erste}.

In the round conformal class we can compare an entire first cluster, rather than only the first gap.
Here we crucially use the twistor-spinor identities of~\cite{Friedrich1990Twistors}.

\begin{theorem}\label{thm:round-conformal-intro}
Let~$g_0$ be the unit round metric on~$\sph^n$,~$n\ge3$, let~$g\in[g_0]$, and put~$N=2^{\lfloor n/2\rfloor}$.
If~$0\le R_g\le n(n-1)$, then
\begin{align}\label{eq:round-first-intro}
 |\lambda_k(g)|\le\frac n2=|\lambda_k(g_0)|,\qquad 0<|k|\le N .
\end{align}
When~$n\not\equiv3\pmod4$, the same conclusion holds under the upper bound~$R_g\le n(n-1)$ alone.
In either case equality for some~$\lambda_k$ in this range holds if and only if~$(\sph^n,g)$ is isometric to the unit round sphere.
\end{theorem}
Theorem~\ref{thm:round-conformal-intro} follows from Theorem~\ref{thm:round-conformal-first-level} and Corollary~\ref{cor:round-upper-symmetry}.
Under~$R_g\le n(n-1)$ alone, the two transported spaces of round Killing spinors show that at least~$2N$ eigenvalues of~$\D_g$ lie in~$[-n/2,n/2]$.
The kernel is zero, so when~$n\not\equiv3\pmod4$, spectral symmetry pairs the eigenvalues in this interval and at least~$N$ have each sign~\cite[Theorem~1.3.7(iv)]{Ginoux2009DiracSpectrum}.
Dimension~$4\ell+3$ enters here for a different reason than in Theorem~\ref{thm:odd-intro}: the dimension no longer forces spectral symmetry, so the total count alone does not ensure that at least~$N$ eigenvalues have each sign.
Under the additional assumption~$R_g\ge0$, we count the two signs separately. 
We don't know whether the scalar upper bound alone suffices for these separate estimates in dimension~$4\ell+3$.
For every~$n\ge3$, the same high-index conclusion~\eqref{eq:round-tail-intro} holds under~$g\in[g_0]$ and~$R_g\le n(n-1)$ unless~$g$ is round; see Proposition~\ref{prop:round-high-indices}.

\medskip 

In Part~\ref{part:conformal} we fix the conformal class and spin structure.
Section~\ref{sec:negative} treats the negative scalar-curvature comparison and its counterpart for the scalar Laplacian above an explicit threshold.
Section~\ref{sec:positive} develops the holomorphic ladders on the sphere and on surfaces of higher genus, and Section~\ref{sec:round-conformal} treats the first cluster and the higher eigenvalues in the round conformal class.
In Section~\ref{sec:flows} we apply signed conformal comparison along the Yamabe flow and, in dimension two, the Ricci flow.

In Part~\ref{part:unrestricted} we allow the conformal class to vary.
Section~\ref{sec:varying-conformal-class} exhibits a jump in the Dirac kernel along a family of hyperbolic metrics with a fixed spin structure and discusses the Dirac eigenvalues on the two-torus.
In Section~\ref{sec:higher-spheres} we use a weighted, locally free circle action on odd-dimensional spheres and prove the first-gap bound for even-dimensional spheres under~$0<\sec_g\le1$.
In even dimensions the same injectivity-radius estimate also implies the requisite volume bound, so the high-index comparison extends to metrics outside the round conformal class.

\subsection*{Conventions}
All manifolds are smooth, closed, and connected unless specified otherwise, and are spin whenever a Dirac operator is used.
Hermitian products are complex linear in the first argument.
All spectral multiplicities and dimensions of spaces of sections are taken over~$\C$.
For the scalar Laplacian these multiplicities also equal the dimensions of the corresponding real-valued eigenspaces.
The notation~$h^q$ for holomorphic cohomology always denotes complex dimension.
We denote the complex spinor bundle over~$M$ by~$\Sigma M$, and write \(\langle\cdot,\cdot\rangle_{\Sigma_gM}\) and \(|\cdot|_{\Sigma_gM}\) for the Hermitian inner product and norm on the spinor bundle associated with the Riemannian metric \(g\).
Clifford multiplication is denoted by~$\gamma$, satisfying the Clifford relation~$\gamma(X)^2=-|X|^2\id$, and the spin connection by~$\snabla$.

The scalar Laplacian is nonnegative:~$\Delta=-\diverg\grad$.
We write~$\sec_g$ for sectional curvature and~$R_g$ for scalar curvature.
In dimension two,~$\sec_g=K_g$ is the Gaussian curvature.
The canonical line bundle is denoted by~$\mathcal K$, and any square root comes with its chosen spin structure.
All coefficients are smooth, so elliptic regularity applies throughout.

We write~$\mu_j(g)$ for the eigenvalues of the scalar Laplacian~$\Delta_g$ and reserve~$\lambda_k(g)$ for those of the Dirac operator.
We distinguish~$\Delta_g$ from the connection Laplacian on a general vector bundle and from the spin Laplacian~$\snabla^*\snabla$.

{\bf Acknowledgement.} We learned the holomorphic ladder method for eigenvalue comparison on two-spheres in an inspiring talk given by Prof. Jianquan Ge, to whom we would like to express our sincere gratitude.

\part{Comparison within a fixed conformal class}\label{part:conformal}

\section{Comparison in negative scalar curvature}\label{sec:negative}

Throughout this section the conformal class and the topological spin structure are fixed.
As in the introduction,~$\lambda_k$,~$k\in\bZ\setminus\{0\}$, denotes a Dirac eigenvalue counted with its complex multiplicity.
No spectral symmetry is assumed, so the statements also apply in dimensions~$4\ell+3$.

\subsection{The conformal variational formula}

Fix a reference metric~$h$ on~$M^n$, write~$g=e^{2u}h$ with~$u\in C^\infty(M,\R)$, and let~$\beta:\Sigma_hM\longrightarrow\Sigma_gM$ be the canonical fiberwise isometry of spinor bundles.
By the conformal covariance of Dirac operators~\cite{Hitchin1974harmonic},
\begin{align}\label{eq:dirac-conformal}
 \D_g\bigl(e^{-(n-1)u/2}\beta\psi\bigr)=e^{-(n+1)u/2}\beta(\D_h\psi).
\end{align}
Consequently, the~$g$-eigenvalue equation is equivalent to the equation
\begin{align}\label{eq:weighted-general}
 \D_h\psi=\lambda e^u\psi,
 \quad\text{on }\Sigma_hM.
\end{align}
The same formula identifies the harmonic-spinor spaces throughout the conformal class.

We recall the reciprocal variational formula from~\cite[Proposition~2.4]{QiuWu2025}, written on the quotient by the kernel so that the trial space is independent of~$u$.
On the complex quotient space
\begin{align}
 \mathscr X=H^1(\Sigma_hM)/\Ker\D_h
\end{align}
set
\begin{align}
 q(\psi)=\int_M\operatorname{Re}\langle\D_h\psi,\psi\rangle_{\Sigma_h}\,\dv_h,
 \qquad
 N_u(\psi)=\int_Me^{-u}|\D_h\psi|_{\Sigma_h}^2\,\dv_h.
\end{align}
Both expressions are well-defined modulo harmonic spinors, and~$N_u$ is positive on every nonzero class.

\begin{lemma}\label{lemma:reciprocal}
For every~$k\ge1$,
\begin{align}
 \frac1{\lambda_k(g)}
 &=\sup_{\substack{E\subset\mathscr X\\\dim_\C E=k}}\ \inf_{0\ne\psi\in E}\frac{q(\psi)}{N_u(\psi)},
 \label{eq:reciprocal-positive}\\
 \frac1{\lambda_{-k}(g)}
 &=\inf_{\substack{E\subset\mathscr X\\\dim_\C E=k}}\ \sup_{0\ne\psi\in E}\frac{q(\psi)}{N_u(\psi)}.
 \label{eq:reciprocal-negative}
\end{align}
The extrema range over all complex~$k$-dimensional subspaces of~$\mathscr X$; no restriction~$q(\psi)>0$ is imposed.
\end{lemma}

We next apply scalar-weight comparison to conformal metrics.

\begin{prop}\label{prop:weight-order}
Let~$g_i=e^{2u_i}h$,~$i=1,2$.
If~$u_1\ge u_2$, then for each $k\in\mathbb Z_*$, we have
\begin{align}
|\lambda_k(g_1)|\le|\lambda_k(g_2)|.
\end{align}
If~$u_1\not\equiv u_2$, the inequality is strict for every~$k\in\mathbb Z_*$.
\end{prop}

\begin{proof}
By~\eqref{eq:weighted-general}, the Dirac eigenvalues of~$g_i$ are the weighted eigenvalues of~$\D_h$ with the positive scalar weights~$a_i=e^{u_i}$.
Since~$a_1\ge a_2$, the asserted ordering follows from~\cite[Theorem~1.3]{QiuWu2025}.
If equality holds for some~$\lambda_k$, then~$a_1=a_2$ by~\cite[Corollary~2.7]{QiuWu2025}, and hence~$u_1=u_2$.
\end{proof}

\begin{cor}\label{cor:conformal-extrema}
For~$g=e^{2u}h$, put~$u_*=\min_Mu$ and~$u^*=\max_Mu$.
Then for each $k\in\mathbb Z_*$,
\begin{align}\label{eq:conformal-positive}
 e^{-u^*}|\lambda_k(h)|\le|\lambda_k(g)|\le e^{-u_*}|\lambda_k(h)|.
\end{align}
If~$u$ is nonconstant, both inequalities are strict for every~$k\ne0$.
\end{cor}

\begin{proof}
Apply Proposition~\ref{prop:weight-order} to~$e^{2u_*}h$,~$g$, and~$e^{2u^*}h$, and use~$\lambda_k(c^2h)=c^{-1}\lambda_k(h)$.
If~$u$ is nonconstant, neither~$u-u_*$ nor~$u^*-u$ vanishes identically, so the strict part applies to both comparisons.
\end{proof}

The analogous estimate for the squared spectrum is~\cite[Proposition~5.2]{Ammann2001spectral}.
Its first nonzero eigenvalue is~$\min\{\lambda_1^2,\lambda_{-1}^2\}$. 
The signed characterization permits us to compare the two signs separately even when the spectrum is not symmetric.

\subsection{Rigidity via maximum principle}

Strict negativity of the reference scalar curvature converts curvature ordering into pointwise ordering of the conformal factor, thus yielding the comparison for all Dirac eigenvalues.

\begin{thm}\label{thm:negative-reference}
Let~$(M^n,h)$ be a closed spin manifold,~$n\ge2$, with negative scalar curvature~$R_h<0$ and let~$g\in[h]$.
For every~$k\ne0$,
\begin{align}
 R_g\ge R_h&\quad\Longrightarrow\quad |\lambda_k(g)|\le|\lambda_k(h)|;\\
 R_g\le R_h&\quad\Longrightarrow\quad |\lambda_k(g)|\ge|\lambda_k(h)|.
\end{align}
In either case, equality for some~$\lambda_k$ holds if and only if~$g=h$.
\end{thm}

The reference scalar curvature~$R_h$ need not be constant, and in the first implication~$R_g$ may change sign.
Neither assertion requires a sectional-curvature bound.

\begin{proof}[Proof of Theorem~\ref{thm:negative-reference}]
Write~$g=e^{2u}h$ for some~$u\in C^\infty(M)$.
The scalar curvatures of~$g$ and~$h$ are related by (see e.g.~\cite{Besse1987,LeeParker1987})
\begin{align}\label{eq:scalar-conformal}
 e^{2u}R_g=R_h+2(n-1)\Delta_hu-(n-1)(n-2)|du|_h^2,
\end{align}
which, for~$n=2$, reduces to
\begin{align}
 e^{2u}K_g=K_h+\Delta_hu.
\end{align}

Suppose first that~$R_g\ge R_h$.
At a minimum point~$x_*$ of~$u$, we have~$du(x_*)=0$ and~$\Delta_hu(x_*)\le0$.
Thus~\eqref{eq:scalar-conformal} implies
\begin{align}
 R_g(x_*)\le e^{-2u_*}R_h(x_*).
\end{align}
If~$u_*<0$, then~$e^{-2u_*}>1$ and~$R_h(x_*)<0$, so
\begin{align}
 R_g(x_*)\le e^{-2u_*}R_h(x_*)<R_h(x_*),
\end{align}
a contradiction to~$R_g\ge R_h$. 
Hence~$u\ge0$, and the eigenvalue comparison follows from Proposition~\ref{prop:weight-order}.

Suppose instead that~$R_g\le R_h$.
At a maximum point~$x^*$, we have~$du(x^*)=0$ and~$\Delta_hu(x^*)\ge0$, whence
\begin{align}
 R_g(x^*)\ge e^{-2u^*}R_h(x^*).
\end{align}
If~$u^*>0$, then~$e^{-2u^*}<1$ and the right-hand side is strictly greater than~$R_h(x^*)$, again a contradiction.
Therefore~$u\le0$, and the reversed comparison follows from Proposition~\ref{prop:weight-order}.

If~$u\not\equiv0$, the strict part of Proposition~\ref{prop:weight-order} applies in either case.
Thus equality for some~$\lambda_k$ forces~$u\equiv0$, that is,~$g=h$.
The converse is immediate.
\end{proof}

\begin{cor}\label{cor:negative-pinching}
Suppose that~$R_h=-c$, where~$c>0$, and~$g\in[h]$.
For every~$k\in\mathbb Z_*$,
\begin{align}
 R_g\ge-cb^2&\quad\Longrightarrow\quad |\lambda_k(g)|\le b|\lambda_k(h)|,\qquad b>0;\\
 R_g\le-ca^2&\quad\Longrightarrow\quad |\lambda_k(g)|\ge a|\lambda_k(h)|,\qquad a>0.
\end{align}
Equality for some~$\lambda_k$ holds precisely when~$g=b^{-2}h$ or~$g=a^{-2}h$, respectively.
In particular,~$-cb^2\le R_g\le-ca^2$, with~$0<a\le b$, implies
\begin{align}\label{eq:negative-two-sided}
 a|\lambda_k(h)|\le|\lambda_k(g)|\le b|\lambda_k(h)|
 \qquad\text{for each }k\in\mathbb Z_*.
\end{align}
\end{cor}

\begin{proof}
For the rescaled reference metric~$h_b=b^{-2}h$, we have
\begin{align}
 R_{h_b}=b^2R_h=-cb^2,\qquad \lambda_k(h_b)=b\lambda_k(h).
\end{align}
Since~$g\in[h_b]$, Theorem~\ref{thm:negative-reference} implies the first comparison and its rigidity.
Apply the other implication to~$h_a=a^{-2}h$ to obtain the second comparison and its equality case.
\end{proof}

On a surface of genus greater than one, fixing the complex structure determines its hyperbolic metric~$g_0$, with~$K_{g_0}=-1$.
For~$g\in[g_0]$ and~$k\ne0$,
\begin{align}
 K_g\ge-1&\quad\Longrightarrow\quad |\lambda_k(g)|\le|\lambda_k(g_0)|;\\
 K_g\le-1&\quad\Longrightarrow\quad |\lambda_k(g)|\ge|\lambda_k(g_0)|.
\end{align}
Equality for some~$\lambda_k$ forces~$g=g_0$.
In higher dimensions, a hyperbolic reference means~$\sec_{g_0}\equiv-1$, and hence~$R_{g_0}=-n(n-1)$. 
The reference spectrum depends on the manifold and its fixed spin structure.

\begin{rmk}\label{rmk:positive-sign}
For~$R_h>0$, a negative minimum or a positive maximum of~$u$ no longer contradicts the curvature ordering.
For example, if~$F$ is a conformal diffeomorphism of the round sphere which is not an isometry, then~$F^*g_0=e^{2u}g_0$ with~$u$ nonconstant, although~$R_{F^*g_0}=R_{g_0}$.
Hence scalar-curvature equality cannot force~$u=0$ in the positive case.
We will use a different argument to prove the sphere comparison in the next section.
\end{rmk}

\begin{rmk}\label{rmk:torus}
On a two-torus, Gauss--Bonnet implies that a Gaussian curvature of definite sign vanishes identically.
In a fixed conformal class the flat metric is unique up to scaling, so it is determined once we prescribe the area.
Without this normalization, the flat metrics~$h$ and~$c^2h$ have the same curvature but~$\lambda_k(c^2h)=c^{-1}\lambda_k(h)$.

Indeed, for metrics with fixed area, Karpukhin--Métras--Polterovich can compare the first positive Dirac eigenvalues of flat and nonflat metrics in the same conformal class, and obtain the optimal metric with minimal first positive eigenvalue.
\begin{prop}[{\cite[Theorem~3.5]{KMP2024}}]\label{prop:rectangular-tori}
Let~$h_b=dx^2+dy^2$ on~$\tor^2_b=\R^2/(\bZ(1,0)+\bZ(0,b))$.
For every smooth~$g\in[h_b]$ with~$\operatorname{Area}(\tor^2_b,g)=b$, we have
\begin{align}\label{eq:rectangular-torus-comparison}
 \lambda_1(g)\geq\lambda_1(h_b)
\end{align}
in either of the following cases:
\begin{enumerate}
\item the spin structure is periodic along both generators and~$b>2\pi$, in which case~$\lambda_1(h_b)=2\pi/b$;
\item the spin structure is periodic along~$(1,0)$ and antiperiodic along~$(0,b)$, and~$b>\pi/2$, in which case~$\lambda_1(h_b)=\pi/b$.
\end{enumerate}
Equality holds if and only if~$g=h_b$.
\end{prop}

\end{rmk}

\subsection{A partial comparison result for the Laplacian in higher dimensions within a conformal class}\label{subsec:laplacian}

As a byproduct of the method, we now compare the scalar Laplacian within a fixed conformal class of negative scalar curvature in dimensions $n\ge3$.
The surface comparisons are treated together with the Dirac operator in Section~\ref{sec:positive}, while the spherical examples and large-index estimates appear in Section~\ref{sec:higher-spheres}.
Write
\begin{align}
0=\mu_0(g)<\mu_1(g)\le\mu_2(g)\le\cdots
\end{align}
for the eigenvalues of~$\Delta_g=-\operatorname{div}_g\nabla_g$ on complex-valued functions, counted with complex multiplicities.
These multiplicities agree with the real dimensions of the real-valued eigenspaces.
The statements below do not require a spin structure.

For~$n\ge3$, the scalar Rayleigh quotient under~$g=e^{2u}g_0$ is
\begin{align}\label{eq:lap-rayleigh-n}
\frac{\int e^{(n-2)u}|df|_{g_0}^2\,\dv_{g_0}}{\int e^{nu}|f|^2\,\dv_{g_0}}.
\end{align}
Unlike the two-dimensional quotient, both weights vary.
Above an explicit spectral threshold, however, we can still compare the two metrics by using the conformal Laplacian.
The conformal Laplacian and its covariance are classical; see~\cite{LeeParker1987}.
We include the deduction in full.

\begin{thm}\label{thm:lap-threshold}
Let~$(M^n,g_0)$ be closed and connected,~$n\ge3$, with constant scalar curvature~$R_{g_0}=-n(n-1)$, and let~$g\in[g_0]$.
Put
\begin{align}
c_n=\frac{n(n-2)}4.
\end{align}
For every index~$j$ with~$\mu_j(g_0)\ge c_n$,
\begin{align}
R_g\ge R_{g_0}&\quad\Longrightarrow\quad\mu_j(g)\le\mu_j(g_0);\\
R_g\le R_{g_0}&\quad\Longrightarrow\quad\mu_j(g)\ge\mu_j(g_0).
\end{align}
In either comparison, equality for some~$\mu_j$ satisfying~$\mu_j(g_0)\ge c_n$ holds if and only if~$g=g_0$.

\end{thm}

In particular, if~$g_0$ is hyperbolic, meaning~$\sec_{g_0}\equiv-1$, the assumptions~$\sec_g\ge-1$ and~$\sec_g\le-1$ imply the respective conclusions.
Only the finitely many eigenvalues below~$c_n$ can violate the stated ordering.

\begin{proof}[Proof of Theorem~\ref{thm:lap-threshold}]
Write~$g=e^{2u}g_0$ and set~$c=(n-2)/(4(n-1))$, so~$cR_{g_0}=-c_n$.
As in the proof of Theorem~\ref{thm:negative-reference}, we have~$u\ge0$ when~$R_g\ge R_{g_0}$ and~$u\le0$ in the reverse case.
For example, at a minimum point~$x_*$ with~$u_*=u(x_*)<0$, where~$\Delta_{g_0}u(x_*)\le0$, the formula
\begin{align}
R_g=e^{-2u}\left(R_{g_0}+2(n-1)\Delta_{g_0}u-(n-1)(n-2)|du|_{g_0}^2\right)
\end{align}
would imply~$R_g(x_*)<R_{g_0}(x_*)$.

Writing~$\phi=e^{-(n-2)u/2}v$, we obtain from the conformal covariance of~$\Delta+cR$ that
\begin{align}\label{eq:lap-transformed}
\frac{\int|d\phi|_g^2\,\dv_g}{\int|\phi|^2\,\dv_g}
=
\frac{\int(|dv|_{g_0}^2-c_n|v|^2)\,\dv_{g_0}}{\int e^{2u}|v|^2\,\dv_{g_0}}
-c\frac{\int R_ge^{2u}|v|^2\,\dv_{g_0}}{\int e^{2u}|v|^2\,\dv_{g_0}}.
\end{align}
Suppose first~$R_g\ge R_{g_0}$, and take~$v$ in the span of an orthonormal list of the first~$j+1$ scalar eigenfunctions of~$g_0$.
Then
\begin{align}
\int(|dv|_{g_0}^2-c_n|v|^2)\dv_{g_0}\le(\mu_j(g_0)-c_n)\int|v|^2\dv_{g_0}
\le(\mu_j(g_0)-c_n)\int e^{2u}|v|^2\dv_{g_0}.
\end{align}
The last inequality uses both~$u\ge0$ and~$\mu_j(g_0)\ge c_n$.
The last term of~\eqref{eq:lap-transformed} is at most~$c_n$.
The invertible multiplication~$v\mapsto\phi$ preserves the dimension of this trial space, so by min--max we have~$\mu_j(g)\le\mu_j(g_0)$.

For the reverse inequality, take the codimension-$j$ space of~$v$ orthogonal to the first~$j$ scalar eigenfunctions of~$g_0$.
On this space,
\begin{align}
\int(|dv|_{g_0}^2-c_n|v|^2)\dv_{g_0}\ge(\mu_j(g_0)-c_n)\int|v|^2\dv_{g_0}
\ge(\mu_j(g_0)-c_n)\int e^{2u}|v|^2\dv_{g_0},
\end{align}
and the last term of~\eqref{eq:lap-transformed} is at least~$c_n$.
The max--min characterization, again preserving codimension under multiplication, proves~$\mu_j(g)\ge\mu_j(g_0)$.

For rigidity, suppose~$g\ne g_0$.
By the strong maximum principle for the conformal scalar-curvature equation,~$u>0$ in the upper case and~$u<0$ in the lower case.
Also~$R_g-R_{g_0}$ cannot vanish identically, since applying the two orderings simultaneously would force~$g=g_0$.
It therefore has the prescribed strict sign on a nonempty open set.

In the upper case, put~$\mu=\mu_j(g_0)$.
If~$\mu>c_n$, the trial-space inequality is strict because~$u_*>0$.
If~$\mu=c_n$, equality in the upper bound for a trial vector would force both
\begin{align}
\int(|dv|_{g_0}^2-c_n|v|^2)\dv_{g_0}=0,
\qquad \int(R_g-R_{g_0})e^{2u}|v|^2\dv_{g_0}=0.
\end{align}
The first equality, within the initial spectral trial space, places~$v$ in the~$c_n$-eigenspace of~$\Delta_{g_0}$.
The second makes it vanish on a nonempty open set, whence, by unique continuation,~$v=0$.
Compactness of the finite-dimensional unit sphere now makes the upper eigenvalue comparison strict.

For the lower case, suppose for contradiction that~$\mu_j(g)=\mu$.
Transport the span of the first~$j+1$ eigenfunctions of~$g$ into the~$v$-variables.
It intersects the codimension-$j$ trial space used above in a nonzero vector.
Its quotient is at most~$\mu$ by the first description and at least~$\mu$ by the second.
For~$\mu>c_n$, the inequality~$u^*<0$ makes the latter inequality strict, a contradiction.
For~$\mu=c_n$, equality forces~$v$ into the~$c_n$-eigenspace of~$\Delta_{g_0}$ and forces
\begin{align}
\int(R_{g_0}-R_g)e^{2u}|v|^2\dv_{g_0}=0.
\end{align}
Unique continuation again leads to a contradiction.
The case~$g=g_0$ is immediate.
\end{proof}

\section{Holomorphic ladders on surfaces}\label{sec:positive}

We now obtain Dirac eigenvalue comparison and rigidity on the two-sphere by a different argument.
By uniformization, an orientation-preserving diffeomorphism identifies every smooth metric with a metric in the round conformal class.
The unique spin structure and the Dirac eigenvalues are preserved under the pullback.
The maximum principle used in the preceding section does not order the conformal factor relative to the round metric.
The relevant structure is instead a sequence of first-order operators between holomorphic line bundles.
Recently, Lin--Wang--Xu proved the lower Laplacian eigenvalue comparison on the sphere in~\cite[Theorem 1.8]{LWX2026}, and Ding--Ge--Li obtained the corresponding upper comparison in~\cite[Theorem 1.3]{DingGeLi2026}.
Their scalar ladder uses the anticanonical powers~$\mathcal K^{-m}$, while here for the Dirac problem we use the half-canonical powers~$\mathcal K^{1/2-m}$.
Passing from functions to half-spinors changes the shift at level~$m$ from~$m(m+1)$ to~$m^2$ and the dimension of the holomorphic kernel from~$2m+1$ to~$2m$.

\subsection{Comparison on the two-dimensional sphere}\label{subsec:sphere}

Write~$g_0$ for the metric of constant Gaussian curvature~$1$ on~$\sph^2$.
For reference, the positive Dirac eigenvalues of the round metric are arranged in the following blocks:
\begin{align}\label{eq:sphere-round-blocks}
 \lambda_j(g_0)=k
 \quad\Longleftrightarrow\quad k(k-1)<j\le k(k+1),
 \qquad k=1,2,\cdots.
\end{align}
By chirality, the eigenvalues~$k$ and~$-k$ each have complex multiplicity~$2k$.
We recover this classical formula~\cite{Sulanke1981Berechnung,Baer1996} below by specializing the ladder to~$K_g\equiv1$.

\begin{thm}\label{thm:sphere}
Let~$g$ be any smooth metric on~$\sph^2$.
For every~$k\in\bZ\setminus\{0\}$,
\begin{align}
 K_g\ge1&\ \Longrightarrow\
 |\lambda_k(g)|\ge|\lambda_k(g_0)|,
 \label{eq:sphere-lower}\\
 K_g\le1&\ \Longrightarrow\
 |\lambda_k(g)|\le|\lambda_k(g_0)|.
 \label{eq:sphere-upper}
\end{align}
Under either curvature hypothesis, equality for some~$\lambda_k$ holds if and only if~$(\sph^2,g)$ is isometric to the unit round sphere.
The upper comparison does not require~$K_g$ to be positive.
Here~$\lambda_{-j}(g)=-\lambda_j(g)$ by chirality.
\end{thm}

\subsubsection{The bundles and the first-order operators}

The metric and orientation determine a complex structure on the sphere, and the resulting Riemann surface is biholomorphic to~$\mathbb{P}^1$.
Let~$\mathcal K=\Lambda^{1,0}T^*\sph^2$.
The unique spin structure determines a holomorphic line bundle~$S$, together with~$S^{\otimes2}\simeq\mathcal K$; thus
\begin{align}
 \mathcal K\simeq\mathcal O(-2),\qquad S\simeq\mathcal O(-1).
\end{align}
We use the convention
\begin{align}
 \Sigma^+\sph^2=S,\qquad \Sigma^-\sph^2=\Lambda^{0,1}T^*\sph^2\otimes S,
 \qquad \D=\sqrt2(\bar\partial_S+\bar\partial_S^*).
\end{align}
These standard identifications are recalled in~\cite{LM1989}.
For~$m\ge0$, set
\begin{align}
 L_m=S\otimes\mathcal K^{-m}
      =\mathcal K^{1/2-m}\simeq\mathcal O(2m-1).
\end{align}
The metric identifies the following Hermitian bundles by a complex-linear isometry:
\begin{align}
 I_m:\Lambda^{0,1}T^*\sph^2\otimes L_m\longrightarrow L_{m+1},
\end{align}
where~$I_m$ raises the~$(0,1)$-covector index by the complex-bilinear extension of the metric, identifying~$\Lambda^{0,1}T^*\sph^2$ with~$T^{1,0}\sph^2=\mathcal K^{-1}$.
Define
\begin{align}
 B_m=\sqrt2 I_m\bar\partial_{L_m}:\Gamma(L_m)\longrightarrow\Gamma(L_{m+1}),
 \qquad
 B_m^*:\Gamma(L_{m+1})\longrightarrow\Gamma(L_m).
\end{align}
Adjoints use the displayed Hermitian bundle metrics and the area form~$\dd A_g$.
Thus, after identifying~$\Sigma^-\sph^2$ with~$L_1$,
\begin{align}\label{eq:sphere-Dirac-block}
 \D=\begin{pmatrix}0&B_0^*\\ B_0&0\end{pmatrix}
 \quad\hbox{on }\Gamma(L_0\oplus L_1).
\end{align}
In particular, on~$\Gamma(S)$,
\begin{align}
 B_0^*B_0=2\bar\partial_S^*\bar\partial_S
 =\D^2|_{\Sigma^+\sph^2}=\snabla^*\snabla+\frac{K_g}{2}.
\end{align}
Thus~$B_0^*B_0$ acts on half-spinors; 
in contrast, the anticanonical ladder starts on scalar functions, and its initial operator is the scalar Laplacian.

\subsubsection{Local formulas and the adjacent curvature identity}

In a local holomorphic coordinate~$z=x+iy$, we have 
\begin{align}
 g=e^{2u}|\dd z|^2,\qquad
 |\dd z|_g^2=|\dd\bar z|_g^2=2e^{-2u},\qquad
 K_g=-4e^{-2u}u_{z\bar z}.
\end{align}
Here~$|\dd z|^2$ in the expression for~$g$ means~$\dd x^2+\dd y^2$.
Choose a holomorphic frame~$s_0$ of~$S$ with~$s_0^{\otimes2}=\dd z$, and put~$s_m=s_0\otimes(\p_z)^{\otimes m}$.
Then
\begin{align}
 h_m:=|s_m|_{L_m}^2=(2e^{-2u})^{1/2-m}.
\end{align}
In a holomorphic frame the connection form of the Chern connection has type~$(1,0)$, and metric compatibility identifies it with~$h_m^{-1}\partial h_m=\partial\log h_m$.
Consequently,
\begin{align}
 \nabla^{L_m}_{\p_z}(a s_m)
 =\bigl(a_z+(2m-1)u_z a\bigr)s_m.
\end{align}
For~$m=0$, this is the positive-chirality part of~$\snabla$.
The Dolbeault operator and its adjoint are
\begin{align}
 \bar\partial_{L_m}(a s_m)
   &=a_{\bar z}\,\dd\bar z\otimes s_m,\label{eq:sphere-dbar}\\
 \bar\partial_{L_m}^*(b\,\dd\bar z\otimes s_m)
   &=-2e^{-2u}\bigl(b_z+(2m-1)u_zb\bigr)s_m.
   \label{eq:sphere-dbaradj}
\end{align}
To verify the second formula, integration by parts for compactly supported coefficients yields
\begin{align}
 \langle\bar\partial(a s_m),b\,\dd\bar z\otimes s_m\rangle_{L^2}
 =2\int h_m a_{\bar z}\overline b\,\dd x\,\dd y
 =-2\int a\,\p_{\bar z}(h_m\overline b)\,\dd x\,\dd y.
\end{align}
Since the section pairing has density~$h_me^{2u}$, the coefficient of the adjoint is~$-2e^{-2u}h_m^{-1}\p_z(h_mb)$, as asserted.
The metric identification is
\begin{align}
 I_m(\dd\bar z\otimes s_m)=2e^{-2u}s_{m+1},
 \qquad
 I_m^{-1}(b s_{m+1})=\tfrac12e^{2u}b\,\dd\bar z\otimes s_m.
\end{align}
Indeed,~$h_{m+1}=\tfrac12e^{2u}h_m$, which also verifies unitarity.
Therefore
\begin{align}\label{eq:sphere-local-B}
 \begin{split}
 B_m(a s_m)&=2\sqrt2e^{-2u}a_{\bar z}s_{m+1};\\
 B_m^*(b s_{m+1})
   &=-\sqrt2\bigl(b_z+(2m+1)u_zb\bigr)s_m.
 \end{split}
\end{align}

\begin{lemma}\label{lemma:sphere-partners}
As operators on sections of~$L_{m+1}$,
\begin{align}\label{eq:sphere-commutator}
 B_mB_m^*-B_{m+1}^*B_{m+1}=(2m+1)K_g.
\end{align}
\end{lemma}
\begin{proof}
Applied to~$b s_{m+1}$, the first composition in~\eqref{eq:sphere-commutator} has coefficient
\begin{align}
 -4e^{-2u}
 \bigl(b_{z\bar z}+(2m+1)u_zb_{\bar z}
                    +(2m+1)u_{z\bar z}b\bigr).
\end{align}
The second has coefficient
\begin{align}
 -4\bigl(\p_z(e^{-2u}b_{\bar z})
                  +(2m+3)u_ze^{-2u}b_{\bar z}\bigr)
 =-4e^{-2u}\bigl(b_{z\bar z}+(2m+1)u_zb_{\bar z}\bigr).
\end{align}
Subtracting and using the formula for~$K_g$ proves the identity.
\end{proof}

The factorization is classical in the theory of magnetic operators on surfaces; compare~\cite{FerapontovVeselov2001}.
We have reproved it for completeness with the bundle maps and normalizations fixed.
Introduce the shifted second-order partners
\begin{align}
 A_m&=B_m^*B_m+m^2,\quad\hbox{on }\Gamma(L_m);\\
 C_m&=B_mB_m^*+m^2,\quad\hbox{on }\Gamma(L_{m+1}).
\end{align}
By Lemma~\ref{lemma:sphere-partners}, they satisfy
\begin{align}\label{eq:sphere-defect}
 C_m=A_{m+1}+(2m+1)(K_g-1).
\end{align}
The following diagram is commutative:
\begin{align}
 \begin{tikzcd}[ampersand replacement=\&,column sep=large,row sep=large]
 \Gamma(L_m)\arrow[r,"B_m"]\arrow[d,"A_m"']\&
 \Gamma(L_{m+1})\arrow[d,"C_m"]\\
 \Gamma(L_m)\arrow[r,"B_m"']\&
 \Gamma(L_{m+1}).
 \end{tikzcd}
\end{align}
Its commutativity means~$B_mA_m=C_mB_m$.
The comparison in~\eqref{eq:sphere-defect} takes place on the \emph{same} bundle~$L_{m+1}$.
In particular, the identity is not a composition of Dolbeault operators in one complex: the metric identifications change the holomorphic bundles between successive steps.

\subsubsection{Holomorphic kernels and multiplicities}

\begin{lemma}\label{lemma:sphere-kernels}
For~$m\ge0$,
\begin{align}
 \dim_{\C}\Ker B_m=2m,\qquad \Ker B_m^*=\{0\}.
\end{align}
In particular,~$\Ker\D_g=\{0\}$ for every metric~$g$ on~$\sph^2$.
\end{lemma}
\begin{proof}
Since~$I_m$ is an isomorphism,~$\Ker B_m=H^0(\mathbb{P}^1,\mathcal O(2m-1))$.
The elementary formula
\begin{align}
 h^0(\mathbb{P}^1,\mathcal O(d))=\max\{d+1,0\}
\end{align}
follows by writing a section in the two affine charts: its coefficient in one chart is a polynomial of degree at most~$d$.
The adjoint kernel is identified with~$H^1(\mathbb{P}^1,\mathcal O(2m-1))$, which vanishes by Serre duality:
\begin{align}
 H^1(L_m)\simeq H^0(\mathcal K\otimes L_m^{-1})^*
       \simeq H^0(\mathcal O(-2m-1))^*=0.
\end{align}
Here~$h^0$ denotes the complex dimension.
These are standard applications of Dolbeault theory and Serre duality; see~\cite{Forster1981}.
More explicitly, if~$B_m^*v=0$, then~$\eta=I_m^{-1}v$ satisfies~$\bar\partial^*\eta=0$.
There are no~$(0,2)$-forms, so~$\eta$ is closed.
Since~$H^1(L_m)=0$, we may write~$\eta=\bar\partial s$ and integrate by parts to obtain~$\|\eta\|^2=\langle s,\bar\partial^*\eta\rangle=0$.
At~$m=0$, both kernels vanish, and~\eqref{eq:sphere-Dirac-block} proves the last assertion.
\end{proof}

For a semibounded self-adjoint elliptic operator~$P$ on a bundle, denote its ordered eigenvalues, including zero modes and starting at index~$1$, by~$\nu_j(P)$, counted with complex multiplicities.
The positive spectra of~$B_m^*B_m$ and~$B_mB_m^*$ agree with multiplicities: on a positive~$\sigma$-eigenspace the maps~$B_m$ and~$\sigma^{-1}B_m^*$ are inverse isomorphisms.
Lemma~\ref{lemma:sphere-kernels} therefore implies
\begin{align}\label{eq:sphere-shift}
 \nu_{j+2m}(A_m)=\nu_j(C_m),\qquad j\ge 1,
\end{align}
while
\begin{align}\label{eq:sphere-bottom}
 \nu_1(A_m)=\cdots=\nu_{2m}(A_m)=m^2,\qquad m\ge 1.
\end{align}
At~$m=0$ the spectra of~$A_0$ and~$C_0$ coincide and are strictly positive.
Moreover,
\begin{align}\label{eq:sphere-half-square}
 \nu_j(A_0)=\lambda^2_j(g).
\end{align}
Indeed, for~$A_0s=\sigma s$, the spinors
\begin{align}
 \frac1{\sqrt2}
 \begin{pmatrix}s\\ \pm\sigma^{-1/2}B_0s\end{pmatrix}
\end{align}
are eigenspinors of~$\D$ with eigenvalues~$\pm\sqrt{\sigma}$.
Conversely, either Dirac eigenvalue equation recovers~$A_0s=\sigma s$ from its first component.
In particular, each sign has the multiplicity of the corresponding eigenvalue of~$A_0$, as claimed in~\eqref{eq:sphere-half-square}.
Squaring the full Dirac operator instead would count both signs together.

\subsubsection{The ordered comparison and rigidity}

The induction step combines spectral pairing with curvature comparison on the common bundle~$L_{m+1}$.
By~\eqref{eq:sphere-defect}--\eqref{eq:sphere-shift} and min--max, for every~$m\ge0$ and~$j\ge1$,
\begin{align}\label{eq:sphere-step}
 \nu_{j+2m}(A_m)=\nu_j(C_m)&\ge\nu_j(A_{m+1}),\qquad\text{when } K_g\ge1;\\
 \nu_{j+2m}(A_m)=\nu_j(C_m)&\le\nu_j(A_{m+1}),\qquad\text{when } K_g\le1.
\end{align}
The equality removes the~$2m$ complex zero modes of~$B_m^*B_m$, while the inequality compares~$C_m$ with~$A_{m+1}$ on~$\Gamma(L_{m+1})$.
For example, if~$K_g\ge1$, the quadratic forms differ by~$(2m+1)\int(K_g-1)|s|^2\dd A_g\ge0$.
We claim that, for every~$k\ge1$,
\begin{align}\label{eq:sphere-iterate}
 \nu_{j+k(k-1)}(A_0)&\ge\nu_j(A_k),\qquad\text{when } K_g\ge1;\\
 \nu_{j+k(k-1)}(A_0)&\le\nu_j(A_k),\qquad\text{when }K_g\le1.
\end{align}
The case~$k=1$ is~\eqref{eq:sphere-step} with~$m=0$.
If the lower comparison holds for~$k$, apply it with index~$j+2k$ and then use~\eqref{eq:sphere-step} with~$m=k$ to obtain
\begin{align}\label{eq:sphere-induction-step}
 \nu_{j+k(k+1)}(A_0)
 \ge\nu_{j+2k}(A_k)
 =\nu_j(C_k)
 \ge\nu_j(A_{k+1}).
\end{align}
This proves the induction step for~$K_g\ge1$.
For~$K_g\le1$, reverse both inequalities in~\eqref{eq:sphere-induction-step}.
The cumulative index shift is~$\sum_{m=0}^{k-1}2m=k(k-1)$.
For~$1\le j\le2k$, the right-hand side of~\eqref{eq:sphere-iterate} is~$k^2$ by~\eqref{eq:sphere-bottom}.
The blocks
\begin{align}
 k(k-1)<\ell\le k(k+1),\qquad\text{for each }k\ge1,
\end{align}
exhaust all positive indices.
When~$K_g\equiv1$, we have~$C_m=A_{m+1}$ at every level, so both inequalities in~\eqref{eq:sphere-iterate} are equalities.
In particular, for the round metric and~$1\le j\le 2k$,
\begin{align}
 \nu_{k(k-1)+j}(A_0)=\nu_j(A_k)=k^2.
\end{align}
The blocks exhaust the positive spectrum, and~\eqref{eq:sphere-half-square} identifies each block with the Dirac eigenvalues~$\pm k$, each of complex multiplicity~$2k$.
This recovers~\eqref{eq:sphere-round-blocks} directly from the ladder.
For a general~$g$ satisfying either curvature bound, take square roots in~\eqref{eq:sphere-iterate} on the same block and use chirality for the negative indices.
This proves both comparisons in Theorem~\ref{thm:sphere}.

\begin{figure}[htbp]
\centering
\resizebox{0.8\textwidth}{!}{%
\begin{tikzpicture}[
 x=1cm,y=1cm,
 every node/.style={font=\normalsize},
 value/.style={draw=laddergray!65,rounded corners=3pt,
   fill=white,minimum width=2.85cm,minimum height=.73cm,inner sep=5pt},
 terminal/.style={value,draw=ladderblue,line width=1pt,fill=ladderblue!6},
 pair/.style={-{Stealth[length=2.1mm]},draw=ladderblue,line width=1.1pt},
 compare/.style={-{Stealth[length=2.1mm]},draw=laddergray,line width=1.1pt},
 smallnote/.style={font=\small,text=laddergray},
 rel/.style={fill=white,inner sep=2pt,font=\large}
]
\foreach \x/\m/\d in {0/0/-1,4/1/1,8/2/3,12/3/5}{
  \node[text=ladderblue] at (\x,1.15)
    {$L_{\m}\simeq\mathcal O(\d)$};
}
\draw[laddergray!20] (-1.45,.72)--(13.45,.72);
\node[value] (a0) at (0,0) {$\nu_{j+6}(A_0)$};
\node[value] (c0) at (4,0) {$\nu_{j+6}(C_0)$};
\node[value] (a1) at (4,-1.75) {$\nu_{j+6}(A_1)$};
\node[value] (c1) at (8,-1.75) {$\nu_{j+4}(C_1)$};
\node[value] (a2) at (8,-3.50) {$\nu_{j+4}(A_2)$};
\node[value] (c2) at (12,-3.50) {$\nu_j(C_2)$};
\node[terminal] (a3) at (12,-5.25) {$\nu_j(A_3)=9$};

\draw[pair] (a0.east)--node[rel]{$=$} (c0.west);
\node[smallnote] at (2,.57) {index shift~$0$};
\draw[pair] (a1.east)--node[rel]{$=$} (c1.west);
\node[smallnote] at (6,-1.18) {index shift~$-2$};
\draw[pair] (a2.east)--node[rel]{$=$} (c2.west);
\node[smallnote] at (10,-2.93) {index shift~$-4$};

\draw[compare] (c0.south)--(a1.north);
\node[rel] at (3.57,-.875) {$\ge$};
\node[smallnote,anchor=west] at (4.27,-.875) {$K_g-1\ge0$};
\draw[compare] (c1.south)--(a2.north);
\node[rel] at (7.57,-2.625) {$\ge$};
\node[smallnote,anchor=west] at (8.27,-2.625) {$3(K_g-1)\ge0$};
\draw[compare] (c2.south)--(a3.north);
\node[rel] at (11.57,-4.375) {$\ge$};
\node[smallnote,anchor=west] at (12.27,-4.375) {$5(K_g-1)\ge0$};

\node[anchor=north,text=ladderblue] at (0,-.48)
  {$=\lambda_{j+6}(g)^2$};
\node[smallnote,anchor=north,align=center] at (12,-5.80)
  {Six complex modes\\
   $\dim_{\C}H^0(\sph^2,L_3)=6$};
\node[anchor=west,align=left,text width=8.5cm,font=\small]
  at (-1.40,-4.85) {
  \textcolor{ladderblue}{Horizontal:} spectral pairing.\\[5pt]
  \textcolor{laddergray}{Vertical:} min--max on the same bundle.
  };
\end{tikzpicture}%
}
\caption{The spherical staircase for~$K_g\ge1$ and~$1\le j\le6$. The successive index shifts remove~$0,2,4$ complex holomorphic modes at the shifted levels~$0,1,4$. The terminal cluster gives~$\lambda^2_{j+6}(g)\ge9$. For~$K_g\le1$, reverse the vertical inequalities; the horizontal equalities are unchanged.}
\label{fig:spherical-staircase}
\end{figure}
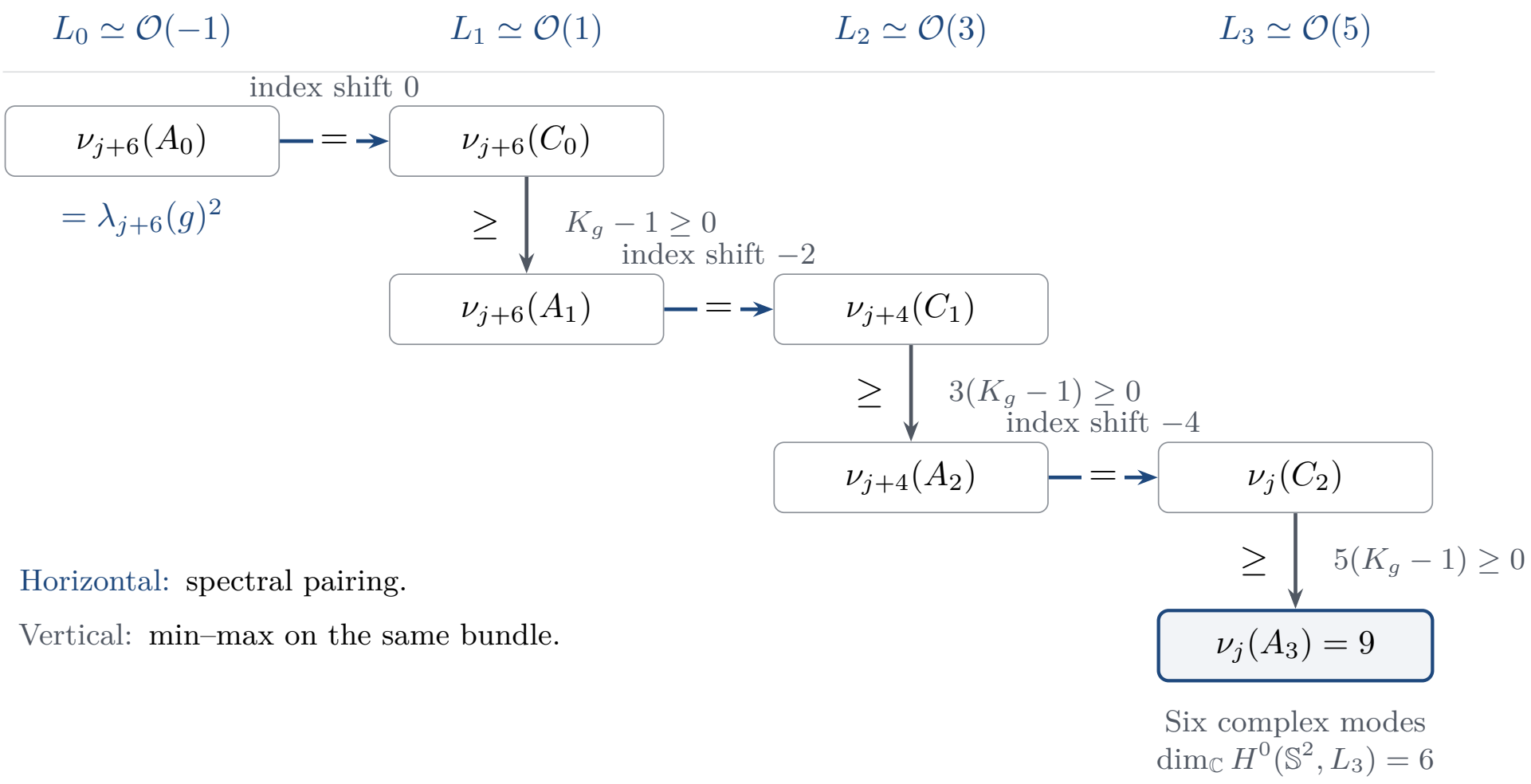

Figure~\ref{fig:spherical-staircase} displays the three-step comparison for the third spectral block.

For the rigidity statement we use the following strict min--max principle; see~\cite[Lemma~3.9]{DingGeLi2026} and its earlier application to the anticanonical ladder in~\cite[Lemma~9.3]{LWX2026}.

\begin{lemma}\label{lemma:strict-potential}
Let~$P$ be a self-adjoint Laplace-type operator on a Hermitian line bundle over a closed connected manifold.
If~$V\ge0$ is smooth and positive somewhere, then~$\nu_j(P+V)>\nu_j(P)$ for every~$j\ge1$.
\end{lemma}

If~$K_g\ge1$ and~$K_g\not\equiv1$, apply this lemma to~\eqref{eq:sphere-defect}.
For each~$m\ge0$ and~$j\ge1$ we then have~$\nu_j(C_m)>\nu_j(A_{m+1})$, so the inequality in the~$K_g\ge1$ case of~\eqref{eq:sphere-step} is strict.
Thus every positive comparison with the round spectrum is strict.
If~$K_g\le1$ and~$K_g\not\equiv1$, apply the same lemma with the nonnegative potential~$(2m+1)(1-K_g)$; the upper comparisons are all strict.
Equality for some~$\lambda_k$ therefore forces~$K_g\equiv1$, and the complete simply connected surface~$(\sph^2,g)$ is isometric to the unit sphere.
This proves Theorem~\ref{thm:sphere}.

\begin{rmk}\label{rmk:scalar-spin-ladder}
For the scalar Laplacian, the ladder acts on the anticanonical powers~$E_m=\mathcal K^{-m}\simeq\mathcal O(2m)$.
Let~$I_m^E:\Lambda^{0,1}T^*\sph^2\otimes E_m\to E_{m+1}$ be the unitary identification obtained by raising the covector index, and define
\begin{align}
 T_m=\sqrt2 I_m^E\bar\partial_{E_m}:\Gamma(E_m)\longrightarrow\Gamma(E_{m+1}).
\end{align}
Since~$E_0$ is trivial,~$T_0^*T_0=2\bar\partial^*\bar\partial=\Delta_g$ on complex-valued functions.
The curvature coefficient is~$2(m+1)$, the complex kernel dimension is~$2m+1$, and the shifted operator is~$T_m^*T_m+m(m+1)$.
The cumulative shift is~$k^2$, yielding the scalar round eigenvalue~$k(k+1)$ with multiplicity~$2k+1$.
For spinors, the curvature coefficient and shifted level become~$2m+1$ and~$m^2$; the complex kernel dimension is~$2m$ and the cumulative index shift is~$k(k-1)$.
Thus it is the change from anticanonical to half-canonical powers that accounts for the different eigenvalues and multiplicities of the two operators.
\end{rmk}

\subsection{Comparison on a hyperbolic surface}\label{subsec:hyperbolic-ladder}

Let~$M$ be a closed Riemann surface of genus~$\gamma>1$, and fix a holomorphic spinor bundle~$S$ with~$S^{\otimes2}\simeq\mathcal K$.
Write~$g_0$ for the unique metric of curvature~$-1$ in a fixed conformal class.
For~$q\in\tfrac12\bZ_{\ge0}$, set~$\mathcal K^q=S^{\otimes2q}$ and define
\begin{align}
 H_{q,g}=2\bar\partial_{\mathcal K^q,g}^{\,*}\bar\partial_{\mathcal K^q},
 \quad\hbox{on }\Gamma(\mathcal K^q).
\end{align}
The holomorphic operator~$\bar\partial_{\mathcal K^q}$ is fixed, while its adjoint and the induced Hermitian bundle metric~$h_{q,g}$ depend on~$g$.
Write~$\mu_j^{(q)}(g)$,~$j\ge1$, for the eigenvalues of~$H_{q,g}$ with zeros included.
We count complex multiplicities at every level, including complex-valued functions when~$q=0$.
Thus~$H_{0,g}=\Delta_g$ and~$\mu_j^{(0)}(g)=\mu_{j-1}(g)$, where~$0=\mu_0(g)<\mu_1(g)\le\cdots$ is the scalar spectrum.
For~$q=\tfrac12$, the operator~$H_{1/2,g}$ is the positive-chirality part of~$\D_g^2$.

\begin{thm}\label{thm:hyperbolic-twists}
For every~$g\in[g_0]$,~$q\in\tfrac12\bZ_{\ge0}$,~$j\ge1$, and~$a,b>0$,
\begin{align}
 K_g\ge-b^2&\ \Longrightarrow\
 \mu_j^{(q)}(g)\le b^2\mu_j^{(q)}(g_0),\label{eq:hyp-upper}\\
 K_g\le-a^2&\ \Longrightarrow\
 \mu_j^{(q)}(g)\ge a^2\mu_j^{(q)}(g_0).\label{eq:hyp-lower}
\end{align}
Equality for a positive eigenvalue occurs precisely for~$g=b^{-2}g_0$ in~\eqref{eq:hyp-upper}, respectively~$g=a^{-2}g_0$ in~\eqref{eq:hyp-lower}.
\end{thm}

The factorization for constant-curvature surfaces and its holomorphic spectral meaning are classical; see~\cite{FerapontovVeselov2001,LopezAlmoroxTejeroPrieto2006}.
We reprove the needed identities for completeness.
The initial scalar and Dirac comparisons also hold relative to a variable negatively curved metric, as the following argument shows.

\subsubsection*{The initial levels}

\begin{prop}\label{prop:surface-variable-reference}
Let~$h$ be a smooth metric with~$K_h<0$, and let~$g\in[h]$.
For~$q=0,\tfrac12,1$ and every~$j\ge1$,
\begin{align}\label{eq:surface-variable-comparison}
 K_g\ge K_h&\ \Longrightarrow\ \mu_j^{(q)}(g)\le\mu_j^{(q)}(h),\\
 K_g\le K_h&\ \Longrightarrow\ \mu_j^{(q)}(g)\ge\mu_j^{(q)}(h).
\end{align}
Equality for a positive eigenvalue holds if and only if~$g=h$.
In particular, these comparisons apply to every positive scalar eigenvalue and every~$|\lambda_k|$, for~$k\ne0$, with the same rigidity.
\end{prop}
\begin{proof}
Write~$g=e^{2u}h$.
With our nonnegative Laplacian convention, the curvature equation is
\begin{align}\label{eq:surface-variable-curvature}
 K_g=e^{-2u}(K_h+\Delta_hu).
\end{align}
Suppose first that~$K_g\ge K_h$.
At a minimum point~$x_*$ of~$u$, we have~$\Delta_hu(x_*)\le0$.
If~$u_*<0$, then by~\eqref{eq:surface-variable-curvature} we would have
\begin{align}
 K_g(x_*)\le e^{-2u_*}K_h(x_*)<K_h(x_*),
\end{align}
a contradiction.
Hence~$u\ge0$.
Moreover,
\begin{align}
 \Delta_hu+c(x)u\ge0,\qquad
 c(x)=-K_h(x)\begin{cases}(e^{2u(x)}-1)/u(x),&u(x)\ne0,\\2,&u(x)=0,\end{cases}
\end{align}
where~$c>0$ is smooth.
The strong minimum principle now implies that either~$u\equiv0$ or~$u>0$ everywhere.
If~$K_g\le K_h$, we argue instead at a maximum point~$x^*$ and conclude that~$u\le0$.
The function~$w=-u\ge0$ satisfies
\begin{align}
 \Delta_hw+\widetilde c(x)w\ge0,\qquad
\widetilde c(x)=-K_h(x)\begin{cases}(1-e^{-2w(x)})/w(x),&w(x)\ne0,\\2,&w(x)=0.\end{cases}
\end{align}
The same minimum principle shows that either~$u\equiv0$ or~$u<0$ everywhere.

The bundle metrics satisfy~$h_{q,g}=e^{-2qu}h_{q,h}$.
Let~$\widehat h_{q,h}$ denote the induced Hermitian metric on~$\Lambda^{0,1}T^*M\otimes\mathcal K^q$.
The Rayleigh quotient of a nonzero section~$s$ is therefore
\begin{align}\label{eq:hyp-rayleigh}
 \mathcal R_{q,g}(s)
 &=\frac{\langle H_{q,g}s,s\rangle_{L^2(h_{q,g})}}{\|s\|_{L^2(h_{q,g})}^2}\\
 &=\frac{2\displaystyle\int_Me^{-2qu}|\bar\partial s|_{\widehat h_{q,h}}^2\,\dd A_h}
 {\displaystyle\int_Me^{2(1-q)u}|s|_{h_{q,h}}^2\,\dd A_h}.
\end{align}
For~$0\le q\le1$, the condition~$u\ge0$ decreases the numerator and increases the denominator, so~$\mathcal R_{q,g}(s)\le\mathcal R_{q,h}(s)$.
The reverse inequality holds if~$u\le0$.
Since the zero modes are fixed,~\eqref{eq:surface-variable-comparison} follows by min--max.
If~$u>0$ everywhere, compactness implies~$u_*>0$ and the sharper estimate
\begin{align}
 \mathcal R_{q,g}(s)\le e^{-2u_*}\mathcal R_{q,h}(s).
\end{align}
Thus every positive eigenvalue is strictly smaller.
If~$u<0$ everywhere, then~$e^{-2u^*}>1$ and~$\mathcal R_{q,g}(s)\ge e^{-2u^*}\mathcal R_{q,h}(s)$, so every positive eigenvalue is strictly larger.

For~$q=\tfrac12$, put~$d_S=h^0(M,S)$.
Then for every $j\ge1$,
\begin{align}
 \dim_{\C}\Ker H_{1/2,g}=d_S,\qquad
 \dim_{\C}\Ker\D_g=2d_S,\qquad
 \mu_{j+d_S}^{(1/2)}(g)=\lambda^2_j(g).
\end{align}
The last identity follows from the block decomposition of~$\D$ exactly as in~\eqref{eq:sphere-half-square}.
By chirality,~$\lambda_{-j}=-\lambda_j$, and taking square roots proves the assertion for all~$|\lambda_k|$.
\end{proof}

Apply Proposition~\ref{prop:surface-variable-reference} with~$h=b^{-2}g_0$ or~$h=a^{-2}g_0$.
Since constant rescaling multiplies the eigenvalues of~$H_q$ by the inverse square scale, this proves Theorem~\ref{thm:hyperbolic-twists} for~$q=0,\tfrac12,1$, including rigidity.
In terms of the conformal factor~$g=e^{2u}g_0$, the conclusions of the maximum principle read
\begin{align}\label{eq:hyp-factor-bounds}
 K_g\ge-b^2\ \Longrightarrow\ u\ge-\log b,\qquad
 K_g\le-a^2\ \Longrightarrow\ u\le-\log a,
\end{align}
with a strict inequality everywhere unless the corresponding homothety holds.
The scalar conclusion is the conformal Rayleigh-quotient comparison; the half-canonical conclusion is a Dolbeault proof of the surface case of the signed conformal comparison.
Both use the same curvature equation and maximum principle.

\subsubsection*{Propagation to all canonical twists}

Let
\begin{align}
 I_{q,g}:\Lambda^{0,1}T^*M\otimes\mathcal K^q
          \longrightarrow\mathcal K^{q-1},\qquad
 A_{q,g}=\sqrt2I_{q,g}\bar\partial_{\mathcal K^q}.
\end{align}
In a holomorphic coordinate with~$g=e^{2v}|\dd z|^2$ and frame~$s_q=(\dd z)^q$, we compute
\begin{align}
 \begin{split}
 A_{q,g}(a s_q)&=2\sqrt2e^{-2v}a_{\bar z}s_{q-1},\\
 A_{q,g}^*(b s_{q-1})
    &=-\sqrt2\bigl(b_z+2(1-q)v_zb\bigr)s_q.
 \end{split}
\end{align}
The calculation in Lemma~\ref{lemma:sphere-partners}, now with exponent~$q$, yields the identities on the indicated bundles:
\begin{align}\label{eq:hyp-partner}
 \begin{aligned}
 A_{q,g}^*A_{q,g}&=H_{q,g},
       &&\hbox{on }\Gamma(\mathcal K^q);\\
 A_{q,g}A_{q,g}^*&=H_{q-1,g}-2(q-1)K_g,
       &&\hbox{on }\Gamma(\mathcal K^{q-1}).
 \end{aligned}
\end{align}
For~$q>1$, we have by Serre duality
\begin{align}
 H^1(M,\mathcal K^q)
 \simeq H^0(M,\mathcal K^{1-q})^*=0,
\end{align}
because the last bundle has negative degree.
Thus~$A_{q,g}^*$ has zero kernel, and the full spectrum of the second operator in~\eqref{eq:hyp-partner} is exactly the positive spectrum of~$H_{q,g}$.
The number of omitted zero modes is fixed:
\begin{align}
 h_q:=\dim_{\C}H^0(M,\mathcal K^q)=(2q-1)(\gamma-1),\qquad q>1.
\end{align}
For the initial levels~$q=0,\tfrac12,1$, the complex kernel dimensions are~$1,h^0(M,S),\gamma$, respectively.
Here~$h^0$ denotes complex dimension in the holomorphic Riemann--Roch formula.

By spectral pairing and~\eqref{eq:hyp-partner}, for every~$q>1$ and~$j\ge1$ we have
\begin{align}\label{eq:hyp-spectral-pairing}
 \mu_{j+h_q}^{(q)}(g)
 =\nu_j(A_{q,g}A_{q,g}^*)
 =\nu_j\bigl(H_{q-1,g}-2(q-1)K_g\bigr).
\end{align}
Thus the two curvature assumptions lead to the explicit induction formulas
\begin{align}\label{eq:hyp-induction-step}
 \mu_{j+h_q}^{(q)}(g)=\nu_j(A_{q,g}A_{q,g}^*)
 &\le\mu_j^{(q-1)}(g)+2(q-1)b^2,
 &&K_g\ge-b^2;\\
 \mu_{j+h_q}^{(q)}(g)=\nu_j(A_{q,g}A_{q,g}^*)
 &\ge\mu_j^{(q-1)}(g)+2(q-1)a^2,
 &&K_g\le-a^2.
\end{align}
The inequalities follow from min--max on the same bundle~$\mathcal K^{q-1}$.
For the reference metric,~$K_{g_0}=-1$ makes the curvature term constant, so
\begin{align}\label{eq:hyp-reference-recurrence}
 \mu_{j+h_q}^{(q)}(g_0)=\mu_j^{(q-1)}(g_0)+2(q-1).
\end{align}

Suppose the upper comparison is known at level~$q-1$.
Combining~\eqref{eq:hyp-induction-step} with~\eqref{eq:hyp-reference-recurrence} yields
\begin{align}\label{eq:hyp-upper-induction}
 \mu_{j+h_q}^{(q)}(g)
 &\le\mu_j^{(q-1)}(g)+2(q-1)b^2
 \le b^2\bigl(\mu_j^{(q-1)}(g_0)+2(q-1)\bigr)
 =b^2\mu_{j+h_q}^{(q)}(g_0).
\end{align}
Likewise, if the lower comparison holds at level~$q-1$, then
\begin{align}\label{eq:hyp-lower-induction}
 \mu_{j+h_q}^{(q)}(g)
 &\ge\mu_j^{(q-1)}(g)+2(q-1)a^2
 \ge a^2\bigl(\mu_j^{(q-1)}(g_0)+2(q-1)\bigr)
 =a^2\mu_{j+h_q}^{(q)}(g_0).
\end{align}
The index~$j+h_q$ omits exactly the kernel of~$H_{q,g}$, while~$\mu_j^{(q-1)}$ includes the kernel at level~$q-1$.
Adding back the same~$h_q$ zero modes on both sides therefore proves comparison at every index at level~$q$.
The initial levels~$q=\tfrac12$ and~$q=1$ were proved above, so these recurrences propagate the comparison along the half-integer and integer chains, respectively.
Together with the scalar level~$q=0$, this proves both inequalities in Theorem~\ref{thm:hyperbolic-twists}.

For the rigidity when~$q>1$, if~$K_g\ge-b^2$ and~$K_g\not\equiv-b^2$, Lemma~\ref{lemma:strict-potential} applied to
\begin{align}
 H_{q-1,g}+2(q-1)b^2
  -\bigl(H_{q-1,g}-2(q-1)K_g\bigr)
  =2(q-1)(K_g+b^2)
\end{align}
makes the first inequality in~\eqref{eq:hyp-upper-induction} strict at every index.
The lower case is identical.
If~$K_g$ is the corresponding constant, uniqueness of the constant negative curvature metric in the fixed conformal class implies~$g=b^{-2}g_0$ or~$g=a^{-2}g_0$.
Conversely, the eigenvalues scale by the stated factor under either homothety, which proves the equality characterization in Theorem~\ref{thm:hyperbolic-twists}.

\begin{rmk}
In the case of hyperbolic surfaces, no exact Dirac spectrum is computed or needed.
The fixed complex and spin structures keep the holomorphic kernels fixed, the maximum principle treats the initial scalar and spinor levels, and the partner identities propagate the comparison to higher twists.
For the reference metric, the curvature term in~\eqref{eq:hyp-partner} is constant; for variable curvature it is a multiplication operator.
Thus Proposition~\ref{prop:surface-variable-reference} treats the scalar and Dirac levels for variable~$K_h<0$, while the induction above uses a constant-curvature reference at the higher levels.
Ordering the eigenvalues of~$H_{q-1,g}$ and~$H_{q-1,h}$ alone does not compare their sums with a variable curvature multiplication operator.
\end{rmk}

\begin{figure}[htbp]
\centering
\resizebox{0.8\textwidth}{!}{%
\begin{tikzpicture}[
 x=1cm,y=1cm,
 every node/.style={font=\normalsize},
 value/.style={draw=laddergray!65,rounded corners=3pt,
   fill=white,minimum width=2.85cm,minimum height=1.03cm,
   inner sep=5pt,align=center},
 pair/.style={-{Stealth[length=2.1mm]},draw=ladderblue,line width=1.1pt},
 rise/.style={-{Stealth[length=2.1mm]},draw=laddergray,line width=1.1pt},
 smallnote/.style={font=\small,text=laddergray},
 rel/.style={fill=white,inner sep=2pt,font=\large}
]
\foreach \x/\q in {0/{1/2},4/{-1/2},8/{-3/2},12/{-5/2}}{
  \node[text=ladderblue] at (\x,1.25) {$\mathcal K^{\q}$};
}
\draw[laddergray!20] (-1.45,.85)--(13.45,.85);
\node[value] (a0) at (0,0) {$B_0^*B_0$\\[3pt]$\lambda^2$};
\node[value] (c0) at (4,0) {$B_0B_0^*$\\[3pt]$\lambda^2$};
\node[value] (a1) at (4,-1.95) {$B_1^*B_1$\\[3pt]$\lambda^2+1$};
\node[value] (c1) at (8,-1.95) {$B_1B_1^*$\\[3pt]$\lambda^2+1$};
\node[value] (a2) at (8,-3.90) {$B_2^*B_2$\\[3pt]$\lambda^2+4$};
\node[value] (c2) at (12,-3.90) {$B_2B_2^*$\\[3pt]$\lambda^2+4$};
\node[value,draw=ladderblue,fill=ladderblue!5] (a3)
 at (12,-5.85) {$B_3^*B_3$\\[3pt]$\lambda^2+9$};

\draw[pair] (a0.east)--node[above=3pt]{$B_0$} (c0.west);
\draw[pair] (a1.east)--node[above=3pt]{$B_1$} (c1.west);
\draw[pair] (a2.east)--node[above=3pt]{$B_2$} (c2.west);
\draw[rise] (c0.south)--(a1.north);
\node[rel,anchor=west] at (4.16,-.975) {$+1$};
\draw[rise] (c1.south)--(a2.north);
\node[rel,anchor=west] at (8.16,-2.925) {$+3$};
\draw[rise] (c2.south)--(a3.north);
\node[rel,anchor=west] at (12.16,-4.875) {$+5$};

\draw[pair,dashed] (a3.east)--node[above=3pt]{$B_3$} (14.65,-5.85);
\draw[rise,dashed] (14.85,-6.03)--(14.85,-6.70);
\node[smallnote,anchor=east] at (14.64,-6.36) {$+7$};
\node[text=ladderblue,font=\Large] at (14.99,-7.03) {$\ddots$};
\node[smallnote,anchor=north] at (0,-.65)
  {$B_0^*B_0\phi_0=\lambda^2\phi_0,\quad\lambda>0$};
\node[anchor=west,align=left,text width=8.7cm,font=\small]
  at (-1.40,-5.65) {
  \textcolor{ladderblue}{Horizontal:} transfer the eigensection by~$B_m$.\\[7pt]
  \textcolor{laddergray}{Vertical:} keep the section and add~$2m+1$.
  };
\end{tikzpicture}%
}
\caption{The forward ladder on the hyperbolic reference metric~$K_{g_0}=-1$. Each box displays an operator and the eigenvalue carried by the chain. Starting from a positive Dirac eigenvalue~$\lambda$, the successive levels are~$\lambda^2+m^2$; the chain therefore does not terminate. The diagram propagates the given spectral value rather than determining it.}
\label{fig:hyperbolic-staircase}
\end{figure}
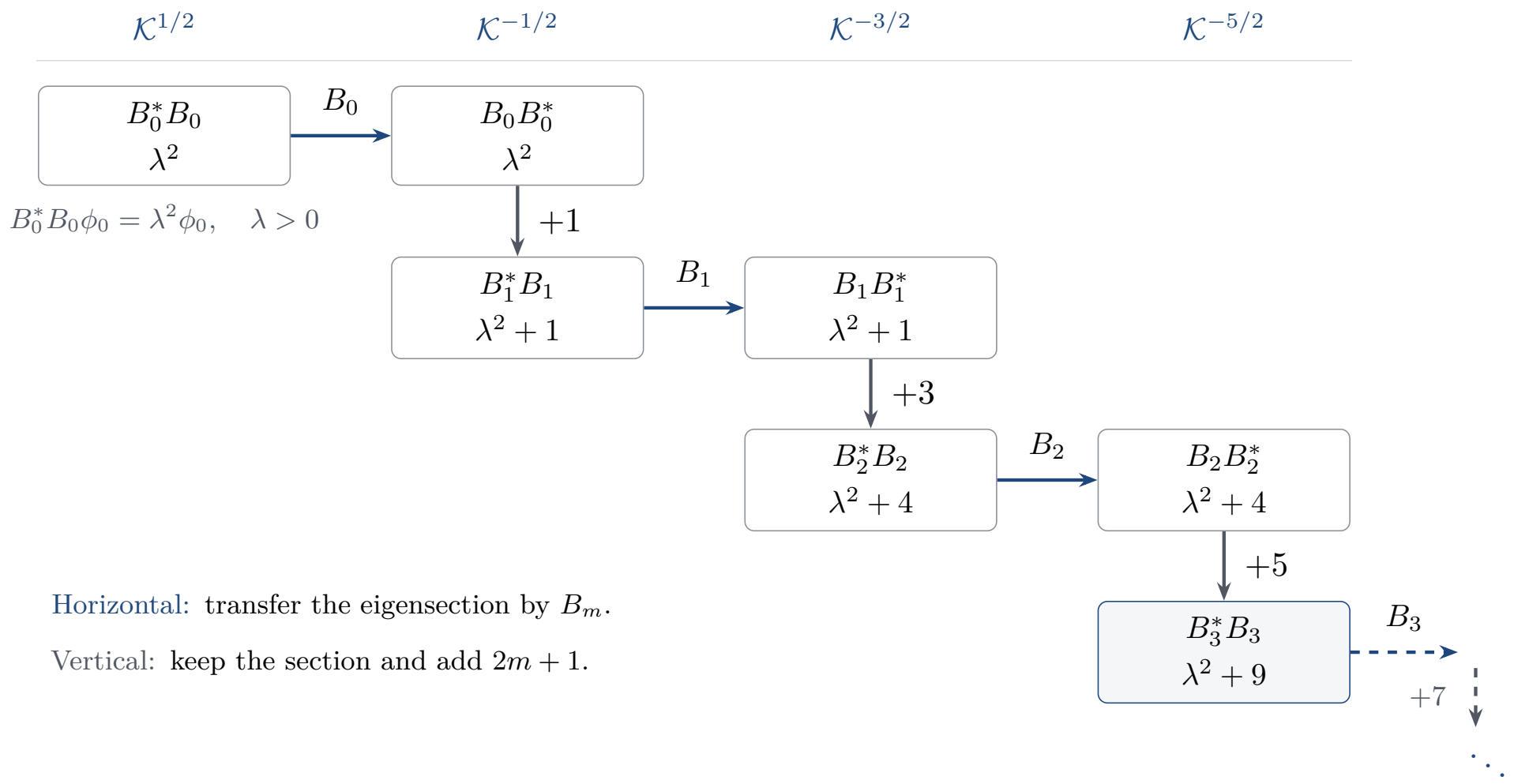

\begin{rmk}\label{rmk:hyperbolic-forward-ladder}
To compare the shape of the hyperbolic construction with Figure~\ref{fig:spherical-staircase}, return to the forward sequence
\begin{align}
 L_m=S\otimes\mathcal K^{-m}=\mathcal K^{1/2-m},\qquad
 B_m=\sqrt2 I_m\bar\partial_{L_m}:
       \Gamma(L_m)\longrightarrow\Gamma(L_{m+1}),
 \qquad m\ge0,
\end{align}
where~$I_m:\Lambda^{0,1}T^*M\otimes L_m\to L_{m+1}$ is the metric-induced unitary identification for~$g_0$.
This is the negative-power direction, whereas the preceding comparison was propagated through the positive canonical powers~$\mathcal K^q$.
Since~$K_{g_0}=-1$, the adjacent identity becomes
\begin{align}
 B_{m+1}^*B_{m+1}=B_mB_m^*+(2m+1).
\end{align}
Take a nonzero eigensection~$\phi_0$ with~$B_0^*B_0\phi_0=\lambda^2\phi_0$, where~$\lambda>0$, and define~$\phi_{m+1}=B_m\phi_m$.
Applying~$B_m$ to the eigenvalue equation and using the adjacent identity, we find inductively that
\begin{align}
 B_m^*B_m\phi_m=(\lambda^2+m^2)\phi_m,\qquad
 \|\phi_{m+1}\|_{L^2(L_{m+1})}^2
   =(\lambda^2+m^2)\|\phi_m\|_{L^2(L_m)}^2>0.
\end{align}
Thus the forward chain in Figure~\ref{fig:hyperbolic-staircase} never reaches a zero mode.
Consistently,~$\deg L_m=(1-2m)(\gamma-1)<0$ for~$m\ge1$, so~$H^0(M,L_m)=0$.
Harmonic spinors are excluded by~$\lambda>0$ and a holomorphic initial section in~$\Ker B_0$ stops at the first step.

These exact eigensection relations concern the constant-curvature reference metric.
For variable curvature the adjacent shift is a multiplication operator, and~$B_m$ need not carry an eigensection to an eigensection of the next unpaired operator.
The infinite chain retains the unknown initial Dirac eigenvalue; it does not yield an explicit formula for the hyperbolic Dirac spectrum.
\end{rmk}

\section{Comparison in the round conformal class}\label{sec:round-conformal}

In this section,~$g_0$ denotes the unit round metric on~$\sph^n$,~$n\ge3$, and~$g=e^{2u}g_0$.
We first compare the lowest round spectral level by transporting the two spaces of round Killing spinors, then use the Yamabe inequality to compare the volume and the higher eigenvalues.

\subsection{The first round spectral level}

Put~$N=2^{\lfloor n/2\rfloor}$.
The round eigenvalues and their complex multiplicities are
\begin{align}\label{eq:round-spectrum-all-dimensions}
 \pm\left(\frac n2+\ell\right),\qquad
 N\binom{n+\ell-1}{\ell},\qquad \ell=0,1,2,\cdots;
\end{align}
see~\cite{Baer1996}.
In particular,~$|\lambda_k(g_0)|=n/2$ for~$0<|k|\le N$.

Let~$\Penr^g$ denote the Penrose (twistor) operator.
For every spinor~$\xi$, the Dirac--Penrose decomposition of its covariant derivative is orthogonal, so pointwise we have
\begin{align}\label{eq:dirac-penrose-norm}
 |\snabla^g\xi|_{T^*M\otimes\Sigma_g}^2
 =|\Penr^g\xi|_{T^*M\otimes\Sigma_g}^2+\frac1n|\D_g\xi|_{\Sigma_g}^2;
\end{align}
see~\cite{Friedrich1990Twistors,Ginoux2009DiracSpectrum}.
This identity holds in every dimension, but a trial space for~$\D_g^2$ counts the two signs together.
To select a sign~$s\in\{1,-1\}$, observe that for~$\lambda\ne0$ we have~$\lambda^2-s(n/2)\lambda\le0$ precisely when~$0<s\lambda\le n/2$.
We therefore apply~$\D_g+s\,n/2$ to the corresponding transported Killing spinors and use the min--max principle for~$\D_g^2-s(n/2)\D_g$.
 This supplies a separate trial space for each sign, including in dimensions~$4\ell+3$.

\begin{thm}\label{thm:round-conformal-first-level}
Let~$g\in[g_0]$ on~$\sph^n$,~$n\ge3$, and suppose that~$0\le R_g\le n(n-1)$.
Then
\begin{align}\label{eq:round-first-level-bound}
 |\lambda_k(g)|\le\frac n2=|\lambda_k(g_0)|,
 \qquad 0<|k|\le N.
\end{align}
Equality for some~$\lambda_k$ in this range holds if and only if~$(\sph^n,g)$ is isometric to the unit round sphere.
\end{thm}

\begin{proof}
Write~$c=n/2$ and~$V=nR_g/[4(n-1)]$, so that~$0\le V\le c^2$.
The conformal covariance formula~\eqref{eq:dirac-conformal} identifies~$\Ker\D_g$ with~$\Ker\D_{g_0}$; the latter is zero by the Schr\"odinger--Lichnerowicz formula.
Consequently,~$\D_g$ is invertible.

Fix~$s\in\{1,-1\}$, and let~$\mathscr K_s$ be the~$N$-dimensional space of round Killing spinors satisfying
\begin{align}\label{eq:round-killing-sign}
 \snabla^{g_0}_X\phi=-\frac s2\gamma_{g_0}(X)\phi,
 \qquad \D_{g_0}\phi=sc\phi.
\end{align}
Every nonzero element of~$\mathscr K_s$ has constant positive length.
Let~$\beta:\Sigma_{g_0}\sph^n\longrightarrow\Sigma_g\sph^n$ be the canonical fiberwise isometry and set
\begin{align}\label{eq:round-transported-twistor}
 \psi=e^{u/2}\beta\phi.
\end{align}
By conformal covariance of the twistor equation~\cite{Friedrich1990Twistors}, we have
\begin{align}\label{eq:round-twistor-equation}
 \Penr^g_X\psi:=\snabla^g_X\psi+\frac1n\gamma_g(X)\D_g\psi=0.
\end{align}
For completeness, combining~$(\Penr^g)^*\Penr^g=(\snabla^g)^*\snabla^g-\D_g^2/n$ with~$\D_g^2=(\snabla^g)^*\snabla^g+R_g/4$ yields
\begin{align}\label{eq:round-twistor-square}
 \D_g^2\psi=\frac{nR_g}{4(n-1)}\psi=V\psi.
\end{align}
Applying~\eqref{eq:dirac-conformal} to~$e^{nu/2}\phi$, we obtain a second identity:
\begin{align}\label{eq:round-twistor-dirac}
 \D_g\psi
 =e^{-u/2}\beta\left(sc\phi+c\gamma_{g_0}(\grad_{g_0}u)\phi\right).
\end{align}
Since Clifford multiplication is skew-Hermitian, this implies the pointwise identities
\begin{align}\label{eq:round-twistor-sign}
 |\psi|_{\Sigma_g}^2=e^u|\phi|_{\Sigma_{g_0}}^2,
 \qquad
 \operatorname{Re}\langle\D_g\psi,\psi\rangle_{\Sigma_g}
 =sc|\phi|_{\Sigma_{g_0}}^2.
\end{align}

Define the modified trial section by
\begin{align}\label{eq:round-modified-trial}
 \eta=(\D_g+sc)\psi.
\end{align}
The map~$\phi\mapsto\eta$ is complex linear and injective, since
\begin{align}
 s\operatorname{Re}\langle\eta,\psi\rangle_{\Sigma_g}
 =c(1+e^u)|\phi|_{\Sigma_{g_0}}^2>0,
 \qquad\phi\ne0.
\end{align}
Its image~$E_s$ therefore has complex dimension~$N$.
By~\eqref{eq:round-twistor-square},~$\D_g\eta=V\psi+sc\D_g\psi$.
Expanding both terms and using~\eqref{eq:round-twistor-sign}, we obtain
\begin{align}\label{eq:round-modified-defect}
 &|\D_g\eta|_{\Sigma_g}^2-sc\operatorname{Re}\langle\D_g\eta,\eta\rangle_{\Sigma_g}\notag\\
 &\quad=(V-c^2)\left(V|\psi|_{\Sigma_g}^2+sc\operatorname{Re}\langle\D_g\psi,\psi\rangle_{\Sigma_g}\right)\notag\\
 &\quad=(V-c^2)(Ve^u+c^2)|\phi|_{\Sigma_{g_0}}^2\le0.
\end{align}
Consider the self-adjoint semibounded operator
\begin{align}\label{eq:round-trial-quotient}
 Q_s=\D_g^2-sc\D_g:\Gamma(\Sigma_g\sph^n)\longrightarrow\Gamma(\Sigma_g\sph^n).
\end{align}
Integrating~\eqref{eq:round-modified-defect} shows that the quadratic form of~$Q_s$ is nonpositive on the~$N$-dimensional complex space~$E_s$.
Hence~$Q_s$ has at least~$N$ nonpositive eigenvalues, counted with complex multiplicities, by the ordinary min--max principle.
On a~$\lambda$-eigenspace of~$\D_g$, the operator~$Q_s$ acts by~$\lambda^2-sc\lambda$.
Since~$\Ker\D_g=0$,
\begin{align}
 \lambda^2-sc\lambda\le0\quad\Longleftrightarrow\quad0<s\lambda\le c.
\end{align}
Thus at least~$N$ eigenvalues of the prescribed sign lie at distance at most~$c$ from zero, which proves~\eqref{eq:round-first-level-bound}.

If~$R_g<n(n-1)$ somewhere, then~$V<c^2$ on a nonempty open set.
The last expression in~\eqref{eq:round-modified-defect} is strictly negative there for every nonzero~$\phi\in\mathscr K_s$, since such a spinor has constant positive length.
The quadratic form of~$Q_s$ is then negative definite on~$E_s$.
Its~$N$th eigenvalue is strictly negative by min--max, so at least~$N$ eigenvalues of~$\D_g$ satisfy~$0<s\lambda<c$.
Consequently, every inequality in~\eqref{eq:round-first-level-bound} is strict.

Equality for some~$\lambda_k$ in the indicated range therefore forces~$R_g\equiv n(n-1)$.
In that case~$V=c^2$, and~\eqref{eq:round-modified-trial} satisfies~$\D_g\eta=sc\eta$ for every~$\eta\in E_s$.
For any such section, the integrated twistor identity reads
\begin{align}
 \int_{\sph^n}|\Penr^g\eta|_{T^*M\otimes\Sigma_g}^2\,\dv_g
 =\frac{n-1}{n}\int_{\sph^n}|\D_g\eta|_{\Sigma_g}^2\,\dv_g
   -\frac14\int_{\sph^n}R_g|\eta|_{\Sigma_g}^2\,\dv_g=0.
\end{align}
Thus~$\snabla^g_X\eta=-s\gamma_g(X)\eta/2$, which is the equality case in Friedrich's estimate~\cite{Friedrich1980erste}.
It follows from the integrability condition of the Killing-spinor equation that~$\Ric_g=(n-1)g$.
As~$g$ is also conformally flat, it has constant sectional curvature~$1$.
Completeness and simple connectedness imply that~$(\sph^n,g)$ is isometric to the unit round sphere.
The converse follows from~\eqref{eq:round-spectrum-all-dimensions}.
\end{proof}

\begin{cor}\label{cor:round-upper-symmetry}
Let~$g\in[g_0]$ on~$\sph^n$,~$n\ge3$, and suppose only that~$R_g\le n(n-1)$.
At least~$2N$ Dirac eigenvalues, counted with complex multiplicities, lie in~$[-n/2,n/2]$; if~$g$ is not isometric to~$g_0$, at least~$2N$ lie in the open interval~$(-n/2,n/2)$.
If~$n\not\equiv3\pmod4$, then
\begin{align}\label{eq:round-first-level-symmetry}
 |\lambda_k(g)|\le\frac n2,\qquad 0<|k|\le N,
\end{align}
and equality for some~$\lambda_k$ in this range holds if and only if~$(\sph^n,g)$ is isometric to the unit round sphere.
\end{cor}

\begin{proof}
Apply~\eqref{eq:round-transported-twistor} to the direct sum~$\mathscr K_+\oplus\mathscr K_-$ of the two round Killing-spinor spaces, and denote its image by~$E$.
Here~$\mathscr K_\pm$ denotes~$\mathscr K_{\pm1}$ in~\eqref{eq:round-killing-sign}.
The transport is an isomorphism onto its image, so~$\dim_\C E=2N$.
Since the twistor equation is linear, every~$\psi\in E$ satisfies~\eqref{eq:round-twistor-square}, and integration by parts yields
\begin{align}\label{eq:round-combined-trial}
 \int_{\sph^n}|\D_g\psi|_{\Sigma_g}^2\,\dv_g
 =\int_{\sph^n}\frac{nR_g}{4(n-1)}|\psi|_{\Sigma_g}^2\,\dv_g
 \le\frac{n^2}{4}\int_{\sph^n}|\psi|_{\Sigma_g}^2\,\dv_g.
\end{align}
By min--max for~$\D_g^2$, we have at least~$2N$ eigenvalues of~$\D_g$ in the stated interval.

If~$R_g<n(n-1)$ on a nonempty open set~$U$, inequality~\eqref{eq:round-combined-trial} is strict for every nonzero~$\psi\in E$.
Indeed, vanishing on~$U$ would mean that its round representative~$\phi_++\phi_-$ vanishes there; applying~$\D_{g_0}$, we also have~$(n/2)(\phi_+-\phi_-)=0$ there.
Both summands would then vanish identically because each has constant length.
The maximum of the Rayleigh quotient on the unit sphere of~$E$ is therefore strictly below~$n^2/4$.
If instead~$R_g\equiv n(n-1)$, the rigidity argument in Theorem~\ref{thm:round-conformal-first-level} identifies~$g$ with the unit round metric up to isometry.

Finally,~$\D_g$ has zero kernel, and its spectrum is symmetric about zero, with complex multiplicities, when~$n\not\equiv3\pmod4$~\cite[Theorem~1.3.7(iv)]{Ginoux2009DiracSpectrum}.
In even dimensions the symmetry follows from chirality; when~$n\equiv1\pmod4$, it follows from the parallel antilinear structure that anticommutes with~$\D_g$.
The eigenvalues in either of the symmetric intervals above therefore occur in pairs of opposite sign.
\end{proof}

Integrating the Dirac--Penrose identity~\eqref{eq:dirac-penrose-norm} and using the Schr\"odinger--Lichnerowicz formula, we obtain for every spinor~$\xi$
\begin{align}
 \|\Penr^g\xi\|_{L^2}^2
 =\frac{n-1}{n}\|\D_g\xi\|_{L^2}^2
 -\frac14\int_{\sph^n}R_g|\xi|_{\Sigma_g}^2\,\dv_g.
\end{align}
For an eigenspinor, the nonnegativity of the left-hand side implies Friedrich's lower bound~$\lambda^2\ge n\min R_g/[4(n-1)]$~\cite{Friedrich1980erste}.
In contrast, the pointwise identity~\eqref{eq:round-twistor-square} holds for twistor spinors, which serve here as trial sections for an upper bound.
In dimensions~$4\ell+3$ the count for~$\D_g^2$ alone does not distinguish the signs.
For the modified trial sections in Theorem~\ref{thm:round-conformal-first-level}, the additional assumption~$R_g\ge0$ makes the second factor in~\eqref{eq:round-modified-defect} positive, and the signs can then be treated separately.

\subsection{Volume comparison and the higher eigenvalues}\label{subsec:sphere-weyl}

Combine the scalar curvature bound with the Weyl asymptotic law we can obtain eigenvalue comparison at high indices, for both Dirac and Laplacian operators. 

\begin{prop}\label{prop:round-high-indices}
Let~$g_0$ be the unit round metric on~$\sph^n$,~$n\ge3$, and let~$g\in[g_0]$ satisfy~$R_g\le n(n-1)$.
Write~$v_g=\vol(\sph^n,g)$ and~$v_0=\vol(\sph^n,g_0)$.
Then~$v_g\ge v_0$, with equality if and only if~$(\sph^n,g)$ is isometric to the unit round sphere.
If~$g$ is not isometric to~$g_0$, then
\begin{align}
 \lim_{k\to\pm\infty}\frac{|\lambda_k(g)|}{|\lambda_k(g_0)|}
 &=\left(\frac{v_0}{v_g}\right)^{1/n}<1,\label{eq:round-high-index-limit}\\
 \lim_{j\to\infty}\frac{\mu_j(g)}{\mu_j(g_0)}
 &=\left(\frac{v_0}{v_g}\right)^{2/n}<1.\label{eq:round-lap-high-index-limit}
\end{align}
For every such nonround metric, there is therefore an integer~$k_0(g)$ such that
\begin{align}\label{eq:round-high-index-bound}
 |\lambda_k(g)|<|\lambda_k(g_0)|,\quad |k|\ge k_0(g);\qquad
 \mu_j(g)<\mu_j(g_0),\quad j\ge k_0(g).
\end{align}
\end{prop}

\begin{proof}
Since~$g\in[g_0]$, we have~$Y(\sph^n,[g])=Y(\sph^n,[g_0])=n(n-1)v_0^{2/n}$~\cite{LeeParker1987}.
Testing its defining quotient in the metric~$g$ with the constant function~$1$, we obtain
\begin{align}\label{eq:round-yamabe-volume}
 n(n-1)v_0^{2/n}
 \le\frac{\displaystyle\int_{\sph^n}R_g\,\dv_g}{v_g^{(n-2)/n}}
 \le n(n-1)v_g^{2/n}.
\end{align}
If the volumes are equal, then~$R_g\equiv n(n-1)$, and the rigidity argument in Theorem~\ref{thm:round-conformal-first-level} implies that~$g$ is isometric to~$g_0$.

We now pass from volume to eigenvalues.
Let~$N_s(T,g)$ count the nonzero eigenvalues of sign~$s$ and absolute value at most~$T$, where~$s\in\{1,-1\}$.
For each nonzero cotangent vector the two eigenspaces of the Dirac principal symbol have the same dimension, namely half the complex spinor rank~$N=2^{\lfloor n/2\rfloor}$.
By the Weyl law for elliptic systems~\cite{Ivrii1982Spectral},
\begin{align}\label{eq:round-signed-weyl}
 N_s(T,g)=\frac{N\omega_n}{2(2\pi)^n}v_gT^n+o(T^n),\qquad\text{as }T\longrightarrow\infty,
\end{align}
where~$\omega_n$ is the volume of the Euclidean unit ball.
Inverting the two counting asymptotics for~$g$ and~$g_0$ proves~\eqref{eq:round-high-index-limit} for each sign, regardless of whether the full Dirac spectrum is symmetric.
For the scalar Laplacian the counting function is asymptotic to~$\omega_n v_g T^{n/2}/(2\pi)^n$~\cite{Chavel1984}, whence~\eqref{eq:round-lap-high-index-limit} follows in the same way.
\end{proof}

\begin{rmk}\label{rmk:weyl-without-conformal-class}
Only the volume estimate in Proposition~\ref{prop:round-high-indices} uses~$g\in[g_0]$.
For any two smooth metrics~$g,h$ on the same closed spin manifold of dimension~$n\ge2$, the signed Dirac and scalar Weyl laws still yield the two limits above, with~$v_0$ replaced by~$\vol(M,h)$ and with~$h$ as the reference metric.
Thus~$\vol(M,g)>\vol(M,h)$ suffices for strict comparison at all sufficiently high indices, however the volume inequality is obtained.
If the volumes are equal, the limits are~$1$ and Weyl's law alone does not order the eigenvalues.
The odd-dimensional metrics of Proposition~\ref{prop:lap-odd-upper} illustrate why the conformal hypothesis matters in the first step: they satisfy~$R_g\le n(n-1)$ but have smaller volume than~$g_0$, so both high-index inequalities are reversed.
\end{rmk}

\begin{cor}\label{cor:round-sectional-first-and-high}
Suppose that~$g\in[g_0]$ on~$\sph^n$,~$n\ge3$, and~$0\le\sec_g\le1$.
Then~$|\lambda_k(g)|\le|\lambda_k(g_0)|$ for~$0<|k|\le N$ and for all sufficiently large~$|k|$.
Unless~$(\sph^n,g)$ is isometric to the unit round sphere, all of these inequalities are strict.
If~$\kappa=\max\sec_g<1$, then
\begin{align}\label{eq:round-sectional-quantitative}
 |\lambda_k(g)|&\le\frac n2\sqrt\kappa,\qquad0<|k|\le N;\notag\\
 \lim_{k\to\pm\infty}\frac{|\lambda_k(g)|}{|\lambda_k(g_0)|}&\le\sqrt\kappa<1 .
\end{align}
\end{cor}

\begin{proof}
Since~$0\le R_g\le n(n-1)$, we can apply Theorem~\ref{thm:round-conformal-first-level} at the first round level and Proposition~\ref{prop:round-high-indices} at high indices.
We have~$\kappa>0$ by~\eqref{eq:round-yamabe-volume}.
Applying the first-level bound to~$\kappa g$ and using~$\lambda_k(\kappa g)=\kappa^{-1/2}\lambda_k(g)$ proves the quantitative estimate for low indices.
For high indices use~$R_g\le n(n-1)\kappa$ in~\eqref{eq:round-yamabe-volume} to obtain~$v_g\ge\kappa^{-n/2}v_0$, and substitute this bound in the volume ratio.
\end{proof}

The threshold~$k_0(g)$ may depend on~$g$.
The preceding corollary compares the first round spectral level and all sufficiently high indices; it does not settle the intervening indices.

\section{Applications in conformal curvature flows}\label{sec:flows}

We fix a topological spin structure on~$M$ and let~$g(t)$ vary in a fixed conformal class.
By Proposition~\ref{prop:weight-order}, when the metric decreases pointwise, the absolute value of every nonzero Dirac eigenvalue is nondecreasing; when the metric increases, the inequalities are reversed.
This is the case for Yamabe flows and the surface Ricci flow. 
Here we start with a general conformal flow.
The comparison above applies at different times of a conformal flow without choosing differentiable eigenvalue branches.

\subsection{Comparison between two times}

Let~$g(t)$ be a smooth family on an interval~$I$ satisfying
\begin{align}\label{eq:general-conformal-flow}
 \p_tg(t)=-2F(t,\cdot)g(t),
\end{align}
where~$F:I\times M\to\R$ is smooth.
Integrating from~$s$ to~$t$, with~$s<t$ in~$I$, we find
\begin{align}
 g(t)=e^{2v_{s,t}}g(s),\qquad
 v_{s,t}(x)=-\int_s^tF(\tau,x)\,\dd\tau.
\end{align}
Conformal covariance identifies the eigenvalue problem at time~$t$ with~$\D_{g(s)}\psi=\lambda e^{v_{s,t}}\psi$ on~$\Sigma_{g(s)}M$.
Define
\begin{align}
 A(s,t)=\int_s^t\min_MF(\tau,\cdot)\,\dd\tau,\qquad
 B(s,t)=\int_s^t\max_MF(\tau,\cdot)\,\dd\tau.
\end{align}

\begin{prop}\label{prop:flow-bounds}
For every~$k\ne0$,
\begin{align}\label{eq:flow-bounds}
 e^{A(s,t)}|\lambda_k(g(s))|
 \le|\lambda_k(g(t))|
 \le e^{B(s,t)}|\lambda_k(g(s))|.
\end{align}
Both inequalities are strict if~$v_{s,t}$ is nonconstant.
\end{prop}

\begin{proof}
Since~$-B(s,t)\le v_{s,t}\le-A(s,t)$, the assertion follows from Corollary~\ref{cor:conformal-extrema}, including its strictness statement.
\end{proof}

\begin{rmk}
For a differentiable normalized eigenpair, the weighted Hellmann--Feynman identity~\cite{QiuWu2026Continuity} takes the form
\begin{align}\label{eq:flow-variation}
 \lambda'(t)=\lambda(t)\int_M F(t,x)|\phi(t,x)|_{\Sigma_{g(t)}}^2\,\dv_{g(t)},\qquad
 \int_M|\phi(t)|_{\Sigma_{g(t)}}^2\,\dv_{g(t)}=1.
\end{align}
Indeed, after fixing a background metric~$h$, the equation becomes~$\D_h\psi=\lambda a(t)\psi$, with~$a_t=-Fa$ and normalization~$\int_M a|\psi|_{\Sigma_h}^2\,\dv_h=1$.
Differentiation and self-adjointness cancel the terms containing~$\psi_t$ and yield~$\lambda'=-\lambda\int_M a_t|\psi|_{\Sigma_h}^2\,\dv_h$.
The conformal identification then implies~\eqref{eq:flow-variation}.
At a multiple eigenvalue, the first derivatives of differentiable branches are determined by the variation form on its finite-dimensional eigenspace.
Proposition~\ref{prop:flow-bounds} applies through crossings without choosing such branches.
\end{rmk}

When the scalar curvature is pinched around a negative constant, we can also control the relative spectral error uniformly in the index.
\begin{cor}\label{cor:flow-pinching}
Let~$R_{g_0}=-c<0$ be constant and~$g(t)\in[g_0]$.
If
\begin{align}
 -c(1+\delta(t))\le R_{g(t)}\le-c(1-\delta(t)),\qquad 0\le\delta(t)<1,
\end{align}
then, for every~$k\ne0$,
\begin{align}\label{eq:flow-rate}
 \sqrt{1-\delta(t)}
 \le\frac{|\lambda_k(g(t))|}{|\lambda_k(g_0)|}
 \le\sqrt{1+\delta(t)}.
\end{align}
In particular,~$\delta(t)\le Ce^{-\rho t}$ implies a relative spectral error of order~$O(e^{-\rho t})$, uniformly in~$k$.
If equality holds for some~$\lambda_k$ in either bound, then~$g(t)$ is the corresponding homothety of~$g_0$; if the two metrics have equal volume, then~$g(t)=g_0$.
\end{cor}

\begin{proof}
Apply Corollary~\ref{cor:negative-pinching} with~$a=\sqrt{1-\delta(t)}$ and~$b=\sqrt{1+\delta(t)}$.
The rate follows from
\begin{align}
 \sqrt{1+\delta}-1=\frac{\delta}{\sqrt{1+\delta}+1},\qquad
 1-\sqrt{1-\delta}=\frac{\delta}{1+\sqrt{1-\delta}}.
\end{align}
The equality statement is the homothety rigidity in that corollary, and equal volume fixes the scale.
\end{proof}

On a surface of genus~$\gamma>1$, normalize the area to~$4\pi(\gamma-1)$.
Hamilton's convergence theorem for normalized Ricci flow~\cite{Hamilton1988} supplies the limiting metric~$g_0$ in the fixed conformal class, with~$K_{g_0}=-1$.
Once the curvature is uniformly close to~$-1$, Corollary~\ref{cor:flow-pinching}, with~$c=2$, turns a bound on~$\|K_{g(t)}+1\|_{L^\infty}$ into relative estimates for the entire nonzero Dirac spectrum.
For normalized Yamabe flow in a negative conformal class, long-time existence and convergence to the constant negative scalar-curvature metric with the prescribed volume are classical~\cite[Theorem~2]{Ye1994}.
The same corollary translates uniform scalar-curvature convergence into quantitative spectral convergence.

\subsection{Monotonicity along Yamabe flow}
\label{subsec:yamabe-flow}

For the unnormalized Yamabe flow, the logarithmic conformal factor satisfies~$\p_tu=-R_{g(t)}/2$.
Its direction of change determines the monotonicity of~$|\lambda_k|$, including in dimensions in which the spectrum need not be symmetric.

\begin{cor}\label{cor:yamabe-all-indices}
Let~$g(t)$ be a smooth Yamabe flow on a closed connected spin manifold of dimension~$n\ge2$, with fixed spin structure.
For the unnormalized equation~$\p_tg=-R_{g(t)}g$, each~$|\lambda_k(g(t))|$, for~$k\ne0$, is nondecreasing if~$R_{g(0)}\ge0$, and nonincreasing if~$R_{g(0)}\le0$.
For the volume-normalized equation~$\p_tg=-(R_{g(t)}-r(t))g$, the same assertions hold for
\begin{align}\label{eq:normalized-yamabe-monotone}
 \exp\!\left(\frac12\int_0^t r(\tau)\,\dd\tau\right)|\lambda_k(g(t))|,
 \qquad
 r(t)=\frac{\int_M R_{g(t)}\,\dv_{g(t)}}{\operatorname{Vol}(M,g(t))}.
\end{align}
Under either initial sign assumption, the corresponding quantities are strictly monotone unless~$g(0)$ is scalar-flat, in which case the flow is stationary.
\end{cor}

\begin{proof}
\emph{Unnormalized flow.}
Write~$R=R_{g(t)}$.
With our nonnegative Laplacian convention the scalar curvature satisfies
\begin{align*}
 \p_t R=-(n-1)\Delta_{g(t)}R+R^2.
\end{align*}
The maximum principle preserves either initial weak sign, and the strong maximum principle makes the sign strict at positive times unless~$R_{g(0)}\equiv0$.
For~$s<t$ we have~$g(t)=e^{2v_{s,t}}g(s)$, where
\begin{align}
 v_{s,t}(x)=-\frac12\int_s^t R_{g(\tau)}(x)\,\dd\tau.
\end{align}
If~$R_{g(0)}\ge0$, then~$v_{s,t}\le0$ and Proposition~\ref{prop:weight-order} implies~$|\lambda_k(g(t))|\ge|\lambda_k(g(s))|$.
For~$R_{g(0)}\le0$ the inequalities reverse.
If the initial scalar curvature is not identically zero, the integral has a strict sign for every~$s<t$ and every~$x$, so both comparisons are strict in their respective cases.

\emph{Volume-normalized flow.}
Here the curvature satisfies
\begin{align*}
 \p_t R=-(n-1)\Delta_{g(t)}R+R\bigl(R-r(t)\bigr),
\end{align*}
so the same maximum-principle argument preserves the initial sign, with the same strictness.
To remove the spatially constant term from the metric equation, set
\begin{align}
 \widehat g(t)=\exp\!\left(-\int_0^t r(\tau)\,\dd\tau\right)g(t).
\end{align}
Since~$\p_t\widehat g=-R_{g(t)}\widehat g$, the logarithmic conformal factor from~$\widehat g(s)$ to~$\widehat g(t)$ is again~$-\tfrac12\int_s^t R_{g(\tau)}\,\dd\tau$.
Applying Proposition~\ref{prop:weight-order} to these two metrics and using the homothety law, we obtain the claimed monotonicity of
\begin{align}
 |\lambda_k(\widehat g(t))|
 =\exp\!\left(\frac12\int_0^t r(\tau)\,\dd\tau\right)|\lambda_k(g(t))|.
\end{align}
As before, the comparison is strict whenever the initial scalar curvature is not identically zero.
A scalar-flat initial metric is stationary for both equations.
\end{proof}

A scalar-flat Yamabe flow is stationary even if the Ricci tensor does not vanish.
For the normalized flow, the derivative~$\p_tu=(r(t)-R_{g(t)})/2$ generally takes both signs, so the argument does not imply monotonicity of the uncorrected~$|\lambda_k|$.
For example, a metric of constant nonzero scalar curvature is stationary under the normalized flow, although the quantity in~\eqref{eq:normalized-yamabe-monotone} is strictly monotone.

\begin{rmk}\label{rem:surface-ricci-flow}
Since~$\operatorname{Ric}_g=K_g g$ and~$R_g=2K_g$ on a surface, the Ricci and Yamabe flows coincide for the corresponding choices of normalization.
For the unnormalized surface Ricci flow~$\p_tg=-2K_g g$, Corollary~\ref{cor:yamabe-all-indices} applies directly to~$|\lambda_k(g(t))|$.

For the area-normalized flow~$\p_tg=2(\overline K-K_g)g$, Gauss--Bonnet and preservation of area imply that
\begin{align}
 r(t)=2\overline K,\qquad
 \overline K=\frac{2\pi\chi(M)}{\operatorname{Area}(M,g(0))}.
\end{align}
The corresponding monotone quantity is thus~$e^{\overline K t}|\lambda_k(g(t))|$.
In each case it is nondecreasing for~$K_{g(0)}\ge0$ and nonincreasing for~$K_{g(0)}\le0$, strictly so unless~$g(0)$ is flat.
The first squared Dirac eigenvalue for surface Ricci flow with nonnegative Gaussian curvature was treated in~\cite{LiYuqiao2022}.
The surface case of Corollary~\ref{cor:yamabe-all-indices} covers all nonzero indices, includes rigidity, and also treats nonpositive curvature.
For the non-strict comparison on surfaces, spectral symmetry allows one to use Ammann's comparison for~$\D^{\,2}$.
Ricci flow in higher dimensions need not remain in a fixed conformal class.
\end{rmk}

\begin{rmk}\label{rmk:equal-area-comparison}
Let~$g$ be a smooth metric on a closed spin surface of genus greater than one, and let~$g_0\in[g]$ have the same area and constant curvature~$K_{g_0}=-\kappa<0$.
Put~$x=(-\min_MK_g)/\kappa\ge1$.
The upper multiplier for squared Dirac eigenvalues in~\cite[Theorem~1.1]{Girao2016} is~$e^{x-1}$, while by Corollary~\ref{cor:negative-pinching} we have
\begin{align}\label{eq:equal-area-upper}
 \frac{|\lambda_k(g)|^2}{|\lambda_k(g_0)|^2}\le x\le e^{x-1},\qquad\text{for each } k\in\mathbb{Z}_*.
\end{align}
This improves the multiplier when~$x>1$ and requires no sign assumption on~$K_g$.
When~$\max_MK_g<0$, we also have
\begin{align}\label{eq:equal-area-lower}
 \frac{-\max_MK_g}{\kappa}\le\frac{|\lambda_k(g)|^2}{|\lambda_k(g_0)|^2}.
\end{align}
The lower estimate in~\cite{Girao2016} also allows sign-changing curvature and is not replaced in that generality by~\eqref{eq:equal-area-lower}.
Equality for some~$\lambda_k$ in either sharp curvature bound forces~$g=g_0$, since the homothety from Corollary~\ref{cor:negative-pinching} must preserve area.
The original flow bounds in~\cite[Theorem~1.1]{Girao2016} are likewise strict at every nonzero index unless~$g=g_0$, as suggested in~\cite[Remark~2.4]{Girao2016}.
Indeed, write~$g_0=e^{2u_\infty}g$ and, in the notation of that paper, let~$r=-2\kappa$ and~$\Lambda_1,\Lambda_2$ be the constants of~\cite[Proposition~2.3]{Girao2016}.
If~$u_\infty$ is nonconstant, integrating the conformal-factor bounds there and applying Corollary~\ref{cor:conformal-extrema}, we find
\begin{align*}
 e^{-\Lambda_1/r}\le\min_M e^{2u_\infty}
 <\frac{|\lambda_k(g)|^2}{|\lambda_k(g_0)|^2}
 <\max_M e^{2u_\infty}\le e^{-\Lambda_2/r}.
\end{align*}
If~$u_\infty$ is constant, equality of the two areas forces~$u_\infty=0$.
\end{rmk}

\part{Comparison without fixing the conformal class}\label{part:unrestricted}

\section{Failure of comparison when the conformal class varies}\label{sec:varying-conformal-class}

\subsection{Surfaces with negative curvature}\label{subsec:surface-negative-curvature}

The assumption of a fixed conformal class in Theorem~\ref{thm:negative-reference} cannot be omitted, even if both metrics have constant Gaussian curvature~$-1$ and the topological spin structure is fixed.
The mechanism is the variation of the dimension of harmonic spinors with the complex structure, established by Hitchin~\cite[Section~2.2, Theorem~2.6]{Hitchin1974harmonic}.
We recall a genus-three example and spell out its spectral consequence for completeness.

\begin{prop}\label{prop:hyperbolic-varying-class}
There are a closed oriented surface~$M$ of genus three, a fixed topological spin structure, and hyperbolic metrics~$g_0,g_\ell$ such that
\begin{align}\label{eq:hyperbolic-kernel-jump}
 K_{g_0}=K_{g_\ell}=-1,\qquad\text{as }g_\ell\longrightarrow g_0\ \text{in }C^\infty,\qquad
 \dim_\C\Ker\D_{g_0}=4,\qquad
 \Ker\D_{g_\ell}=0,
\end{align}
and
\begin{align}\label{eq:hyperbolic-small-eigenvalues}
 0<\lambda_1(g_\ell)=|\lambda_{-1}(g_\ell)|\longrightarrow0
 <\lambda_1(g_0)=|\lambda_{-1}(g_0)|.
\end{align}
Consequently, Gaussian-curvature ordering alone cannot imply either universal ordering of the nonzero Dirac eigenvalues when the conformal class is allowed to vary.
\end{prop}

\begin{proof}
Let~$X_0$ be a hyperelliptic Riemann surface of genus three and let~$\pi:X_0\longrightarrow\mathbb{P}^1$ be its degree-two map.
The hyperelliptic line bundle~$\mathcal L_0=\pi^*\cO(1)$ satisfies
\begin{align}\label{eq:genus-three-theta}
 \mathcal L_0^{\otimes2}\cong\mathcal K_{X_0},\qquad
 \deg\mathcal L_0=2,\qquad
 h^0(X_0,\mathcal L_0)=2.
\end{align}
For instance, these identities follow from~$\pi_*\cO_{X_0}=\cO\oplus\cO(-4)$ and the canonical-bundle formula for a double covering of~$\mathbb{P}^1$ branched at eight points.
Fix the topological spin structure determined by~$\mathcal L_0$.

Choose marked nonhyperelliptic complex structures~$X_\ell$ converging smoothly to~$X_0$.
Such structures exist arbitrarily close to~$X_0$: the hyperelliptic locus has complex dimension~$5$, while the Teichm\"uller space of genus-three surfaces has complex dimension~$6$.
Let~$\mathcal L_\ell$ be the theta characteristic determined by the fixed spin structure on~$X_\ell$.
The parity of~$h^0(X_\ell,\mathcal L_\ell)$ is invariant under deformation~\cite[Theorem~1]{Atiyah1971SpinStructures}; see also~\cite[p.~184]{Mumford1971ThetaCharacteristics}.
Hence it is even.
On a nonhyperelliptic genus-three curve, a degree-two line bundle cannot have two independent holomorphic sections.
Indeed, after removing their common base divisor, the ratio of two such sections would define a nonconstant map to~$\mathbb{P}^1$ of degree at most two.
A map of degree one is impossible for positive genus, and a map of degree two would make the curve hyperelliptic.
It follows that
\begin{align}
 h^0(X_\ell,\mathcal L_\ell)=0.
\end{align}

Let~$g_\ell$ and~$g_0$ be the hyperbolic metrics of these marked conformal structures.
Smooth dependence in uniformization implies that~$g_\ell\to g_0$ in~$C^\infty$.
By the Dolbeault description of the spin Dirac operator and Serre duality, we have
\begin{align}
 \dim_\C\Ker\D_g
 =h^0(X,\mathcal L)+h^1(X,\mathcal L)
 =2h^0(X,\mathcal L),
\end{align}
which proves the kernel assertions in~\eqref{eq:hyperbolic-kernel-jump}.

For completeness, the small eigenvalues can be obtained directly from the ordinary min--max principle for~$\D^2$.
Use the metric identification of spinor bundles to transport the four-dimensional complex space~$\Ker\D_{g_0}$ to a space~$E_\ell\subset\Gamma(\Sigma_{g_\ell}M)$.
Under these identifications the coefficients of the Dirac operators converge smoothly, so on this fixed finite-dimensional trial space we have
\begin{align}\label{eq:kernel-trial-convergence}
 \sup_{0\ne\psi\in E_\ell}
 \frac{\|\D_{g_\ell}\psi\|_{L^2\left(\Sigma_{g_\ell}\right)}^2}
      {\|\psi\|_{L^2\left(\Sigma_{g_\ell}\right)}^2}
 \longrightarrow0.
\end{align}
Because~$\D_{g_\ell}$ is invertible, its smallest squared eigenvalue is positive and is bounded above by the expression in~\eqref{eq:kernel-trial-convergence}.
By chiral symmetry~$\lambda_{-1}(g_\ell)=-\lambda_1(g_\ell)$, which proves~\eqref{eq:hyperbolic-small-eigenvalues}.
The first positive eigenvalue of~$g_0$ is strictly positive because zero is an isolated eigenvalue of finite multiplicity.
Finally, the Gaussian curvatures of~$g_0$ and~$g_\ell$ agree.
Any comparison rule depending only on their ordering would therefore apply after interchanging the metrics, forcing equal first positive eigenvalues in contradiction with~\eqref{eq:hyperbolic-small-eigenvalues}.
\end{proof}

\begin{rmk}\label{rmk:kernel-jump-conformal}
This is a consequence of the classical kernel-jumping construction~\cite[Section~2.2, Propositions~2.4--2.5 and Theorem~2.6]{Hitchin1974harmonic}, together with the deformation invariance of parity~\cite{Atiyah1971SpinStructures,Mumford1971ThetaCharacteristics}.
It also explains why keeping the topological spin structure fixed does not replace the assumption of a fixed conformal class.
The nonzero eigenvalues approaching zero are omitted from the indexed list once they reach the kernel, so the convergence in~\eqref{eq:hyperbolic-small-eigenvalues} does not assert continuity of the first positive eigenvalue across a kernel jump.
Within a fixed conformal class the kernel dimension is constant, and this particular obstruction is absent.
\end{rmk}

\subsection{The two-torus}\label{subsec:two-torus}

Note that any nonflat metric on a two-torus must have a sign-changing curvature. The following proposition thus shows that Theorem~\ref{thm:negative-reference} cannot be relaxed to allowing sign-changing reference curvature (even in a fixed conformal class).

\begin{prop}\label{prop:torus-nonflat}
Let $h$ be any nonflat metric on $\tor^2$. There is a smooth positive function $f$ such that, for all sufficiently small $t>0$, the metrics
$g_t=e^{2tf}h$ satisfy
\begin{align*}
 K_{g_t}<K_h\quad\text{everywhere on }\tor^2;\qquad
 |\lambda_k(g_t)|<|\lambda_k(h)|,\quad\text{for }k\neq0.
\end{align*}
For $g_{-t}$ both inequalities are reversed. The same function $f$ works for every spin structure.
\end{prop}

\begin{proof}
Consider the scalar operator $L_h=\Delta_h-2K_h$. The constant test function has Rayleigh quotient zero by Gauss--Bonnet. Its first eigenvalue $\nu_0$ is strictly negative: otherwise the constant function would minimize the quotient and would satisfy
$L_h1=-2K_h=0$, contrary to nonflatness. Choose its positive ground state $f$, so that $L_hf=\nu_0f$.

The curvature transformation law gives
\[
 K_{g_t}=e^{-2tf}(K_h+t\Delta_hf),\qquad
 \left.\frac{d}{dt}\right|_{t=0}K_{g_t}
 =\Delta_hf-2K_hf=\nu_0f<0.
\]
Compactness makes this inequality uniform, and hence $K_{g_t}<K_h$
for all sufficiently small $t>0$. By conformal comparison,
\[
 |\lambda_k(g_t)|\leq e^{-t\min f}|\lambda_k(h)|
 <|\lambda_k(h)|.
\]
The same argument with $-t$ gives $K_{g_{-t}}>K_h$ and
$|\lambda_k(g_{-t})|\geq e^{t\min f}|\lambda_k(h)|>|\lambda_k(h)|$.
\end{proof}

In the case of varied conformal classes, curvature comparison on a torus does not lead to a Dirac eigenvalue comparison. We write
\begin{align}\label{eq:dirac-eigenvalue-gap}
 \delta(g)=\min\{|\lambda|:\lambda\in\Spec(\D_g)\},
\end{align}
with zero included.
This notation does not assume spectral symmetry.

\begin{prop}\label{prop:torus-espace}
For any nontrivial spin structure on $\tor^2$, there is a family of flat
metrics $h_\varepsilon$ with area one such that
$\delta(h_\varepsilon)\longrightarrow\infty$ as
$\varepsilon\searrow0$. For the trivial spin structure,
$\delta(g)=0$ for every metric $g$.
\end{prop}

\begin{proof}
For a nontrivial spin structure, choose an oriented basis of the lattice such that spinors are antiperiodic in the first generator and periodic in the second. On $\tor^2=\mathbb R^2/\mathbb Z^2$, consider the metric
\[
 h_\varepsilon=\varepsilon^2dx^2+\varepsilon^{-2}dy^2
\]
for~$\eps>0$, which has Gaussian curvature~$K_{h_\varepsilon}=0$ and has unit area. 
From~\cite{Friedrich1984zur} and~\cite[Theorem 2.1.1]{Ginoux2009DiracSpectrum} we know that the Dirac spectrum is
\[
 \operatorname{Spect}(\D_{h_\varepsilon})
 =\left\{\pm2\pi
 \sqrt{\frac{(m+\tfrac12)^2}{\varepsilon^2}
                  +n^2\varepsilon^2}:m,n\in\mathbb Z\right\}.
\]
Consequently,
$\delta(h_\varepsilon)=\pi/\varepsilon\longrightarrow\infty$
as $\varepsilon\searrow0$.

Now fix the trivial spin structure and let $g$ be any smooth metric on $\tor^2$. By uniformization, write $g=e^{2u}h_0$, where $h_0$ is flat.
The trivial spin structure admits a nonzero parallel spinor $\psi$ with respect to~$h_0$, so $\D_{h_0}\psi=0$.
Let $\beta:\Sigma_{h_0}\tor^2\to\Sigma_g\tor^2$ be the canonical fiberwise isometry. Conformal covariance in dimension two gives
\[
 \D_g\bigl(e^{-u/2}\beta\psi\bigr)
 =e^{-3u/2}\beta(\D_{h_0}\psi)=0.
\]
Thus~$e^{-u/2}\beta\psi$ is a nonzero harmonic spinor for $g$ and~$\delta(g)=0$.
\end{proof}

\section{Dirac eigenvalues of higher-dimensional spheres}
\label{sec:higher-spheres}

We now allow the conformal class on~$\sph^n$ to vary.
Under~$0<\sec_g\le1$, the first Dirac eigenvalue of each sign can be arbitrarily far from zero in odd dimensions.
In even dimensions the first gap is at most its round value, and the same upper comparison holds at every sufficiently large absolute index.
The latter conclusion, for both the Dirac operator and the scalar Laplacian, follows from volume comparison.

\subsection{Counterexamples on odd-dimensional spheres}

Let $\delta(g)$ be as in~\eqref{eq:dirac-eigenvalue-gap}.
The following metrics are not required to belong to the round conformal class.

\begin{thm}\label{thm:odd-spheres}
For every odd~$n\ge3$, the sphere~$\sph^n$ admits smooth metrics~$\widehat g_r$,~$0<r\le r_0$, such that
\begin{align}
 0<\sec_{\widehat g_r}\le1,
 \qquad \delta(\widehat g_r)\longrightarrow\infty
 \quad\text{as }r\longrightarrow0.
\end{align}
\end{thm}

The spectral mechanism is due to Ammann--B\"ar~\cite{AmmannBaer1998} and Ammann~\cite{Ammann1998}: under nonprojectable circle collapse, the absence of zero spinorial frequency forces the entire Dirac spectrum away from zero.
Ordinary Hopf collapse already suffices in dimensions~$4\ell+1$; in dimensions~$4\ell+3$, its spin lift is projectable and some eigenvalues remain bounded.
Our additional step is to choose a weighted circle action with odd total weight and positive transverse sectional curvature.
This realizes their mechanism in every odd dimension, including~$\sph^3$, without introducing a new analytic collapse estimate.

\begin{proof}[Proof of Theorem~\ref{thm:odd-spheres}]
Put~$n=2m+1$.
On~$\sph^{2m+1}\subset\C^{m+1}$, let~$g_0$ be the unit round metric and set
\begin{align}
 \eta_0=\sum_{j=0}^{m}(x_j\,\dd y_j-y_j\,\dd x_j),
 \qquad H=\ker\eta_0.
\end{align}
For positive weights~$w=(w_0,\ldots,w_m)$, define
\begin{align}\label{eq:weighted-data}
 f_w=\sum_jw_j|z_j|^2,\qquad
 V_w=\sum_jw_j(x_j\partial_{y_j}-y_j\partial_{x_j}),\qquad
 \theta_w=f_w^{-1}\eta_0 .
\end{align}
Then~$\theta_w(V_w)=1$ and~$\iota_{V_w}\dd\theta_w=0$.
With~$P_wX=X-\theta_w(X)V_w$, the standard weighted Sasakian metric has the form
\begin{align}\label{eq:weighted-metric}
 h_w(X,Y)=f_w^{-1}g_0(P_wX,P_wY),\qquad
 g_{w,r}=h_w+r^2\theta_w^2 .
\end{align}
Here~$h_w$ is a transverse tensor with kernel~$\R V_w$, and~$g_{w,r}$ is a smooth metric on the entire sphere.
The transverse metric is K\"ahler, with
\begin{align}\label{eq:parallel-curvature}
 F_w=\dd\theta_w=2\omega_w,\qquad
 \nabla^{h_w}\omega_w=0;
\end{align}
see~\cite{BoyerGalickiSimanca2008,HeSun2016}.

Choose primitive integer weights
\begin{align}
 w=(N,\ldots,N,N+1),\qquad N\ \text{is even}.
\end{align}
Their sum is odd.
The normalized weights~$b=w/N$ tend to~$(1,\ldots,1)$, so~$h_b+\theta_b^2$ tends smoothly to~$g_0$.
For sufficiently large~$N$, this Sasakian metric has positive sectional curvature, and its transverse metric has positive sectional curvature by the submersion formula.
Since~$h_w=N^{-1}h_b$ transversely, we may fix such an~$N$ and constants
\begin{align}
 0<\kappa_w=\min\sec(h_w)\le\max\sec(h_w)=C_w .
\end{align}
All these choices are fixed while~$r$ tends to zero.

For horizontal orthonormal~$X,Y$ and~$U=r^{-1}V_w$, the sectional curvatures of the connection metric are
\begin{align}\label{eq:sasakian-curvature}
 K_{w,r}(X,Y)=K_{h_w}(X,Y)-3r^2\omega_w(X,Y)^2,\qquad
 K_{w,r}(X,U)=r^2 .
\end{align}
The curvature components with three horizontal arguments and one vertical argument vanish because~$\nabla^{h_w}F_w=0$.
Every two-plane has an orthonormal basis~$X,\cos\alpha\,Y+\sin\alpha\,U$; its sectional curvature is therefore the corresponding convex combination of the two quantities in~\eqref{eq:sasakian-curvature}.
Consequently
\begin{align}\label{eq:weighted-curvature-bounds}
 \min\{\kappa_w-3r^2,r^2\}
 \le\sec(g_{w,r})\le\max\{C_w,r^2\}.
\end{align}

The effective circle action generated by~$V_w$ is
\begin{align}
 e^{it}\cdot(z_0,\ldots,z_m)
       =(e^{iw_0t}z_0,\ldots,e^{iw_mt}z_m).
\end{align}
Each ambient rotation through~$2\pi w_j$ contributes~$(-1)^{w_j}$ to its spin lift.
Thus the lifted time-$2\pi$ map is~$(-1)^{\sum_jw_j}=-1$, also on the intrinsic spin bundle.
The spinorial Lie derivative~$\mathcal L_{V_w}$ consequently has frequencies~$ik$,~$k\in\mathbb Z+\tfrac12$.

We recall the part of the Ammann--B\"ar calculation needed here~\cite[Lemmas~4.3--4.4]{AmmannBaer1998}.
Under the natural identification of spinor bundles,
\begin{align}\label{eq:dirac-splitting}
 \D_{g_{w,r}}
 =\D_H+\frac1r\gamma(U)\mathcal L_{V_w}
       -\frac r4\gamma(U)\gamma(F_w).
\end{align}
The horizontal operator~$\D_H$ is formally self-adjoint, commutes with~$\mathcal L_{V_w}$, and anticommutes with~$\gamma(U)$.
On a frequency space,
\begin{align}
 \left(\D_H+\frac{ik}{r}\gamma(U)\right)^2
       =\D_H^2+\frac{k^2}{r^2}.
\end{align}
Since~$F_w=2\omega_w$, its Clifford action has operator norm at most~$2m$.
Summing over the orthogonal Fourier spaces and estimating the last term in~\eqref{eq:dirac-splitting}, we conclude that
\begin{align}\label{eq:weighted-gap}
 \delta(g_{w,r})\ge\frac1{2r}-\frac m2r.
\end{align}
This is the constant-radius estimate underlying~\cite[Theorem~4.2]{Ammann1998}.
Although that theorem is stated for free actions, the displayed identities hold in the local orbifold quotient charts of our locally free action and define smooth operators on the sphere.
The Fourier decomposition comes from the global spin lift, so the same estimate holds across exceptional orbits.

Choose~$r_0>0$ so that~$3r_0^2<\kappa_w$, and put~$C=\max\{C_w,r_0^2\}$.
Equations~\eqref{eq:weighted-curvature-bounds} and~\eqref{eq:weighted-gap} imply that~$\widehat g_r=Cg_{w,r}$ satisfies~$0<\sec_{\widehat g_r}\le1$ and
\begin{align}
 \delta(\widehat g_r)\ge
 C^{-1/2}\left(\frac1{2r}-\frac m2r\right)\longrightarrow\infty .
\end{align}
\end{proof}

\begin{rmk}
\label{rmk:product-example}
For every~$n\ge3$, equip
\begin{align}
 (M_r,g_r)=\sph^{n-1}(1)\times\sph^1(r)
\end{align}
with the product metric and the antiperiodic spin structure on the circle.
Its Ricci tensor and scalar curvature satisfy
\begin{align}
 \operatorname{Ric}_{g_r}=(n-2)g_{\sph^{n-1}}\oplus0,
 \qquad 0\le\operatorname{Ric}_{g_r}\le(n-2)g_r,
 \qquad R_{g_r}=(n-1)(n-2).
\end{align}
Nevertheless, by the product formula for~$\D^2$ we have
\begin{align}
 \delta(g_r)^2=\frac{(n-1)^2}{4}+\frac1{4r^2}
       \longrightarrow\infty.
\end{align}
Thus even~$0\le\operatorname{Ric}_g\le(n-1)g$ on a general closed spin manifold does not imply the spherical upper comparison~$\delta(g)\le n/2$, or any upper bound depending only on the dimension.
In even dimensions these products also show that Proposition~\ref{prop:even-upper}, as a statement about general closed spin manifolds, fails if~$0<\sec_g\le1$ is replaced by~$0\le\sec_g\le1$.
They are not spheres, and their circle direction has zero Ricci curvature, so they do not decide whether strict positivity is necessary on a sphere or what follows from a positive Ricci lower bound.
\end{rmk}

The scalar Laplacian also violates the upper comparison in odd dimensions, although its first positive eigenvalue need not diverge during collapse.
Write~$0=\mu_0(g)<\mu_1(g)\le\mu_2(g)\le\cdots$ for its eigenvalues on complex-valued functions, counted with complex multiplicities.
The following consequence of classical Hopf collapse and Weyl's law makes the difference precise.

\begin{prop}\label{prop:lap-odd-upper}
Let~$g_0$ be the unit round metric on~$\sph^n$.
For every odd~$n\ge3$, there is a smooth metric~$g$ on~$\sph^n$ with~$0<\sec_g\le1$ such that
\begin{align}
 \mu_j(g)>\mu_j(g_0), \quad&\text{for all sufficiently large }j;\notag\\
 |\lambda_k(g)|>|\lambda_k(g_0)|, \quad&\text{for all sufficiently large }|k|.
\end{align}
\end{prop}

\begin{proof}
Put~$n=2m+1$, and let~$g_r$ be the classical Hopf connection metric on~$\sph^{2m+1}$ with round horizontal metric and fiber radius~$0<r\le1$.
Its sectional curvatures lie between~$r^2$ and~$4-3r^2$~\cite{Besse1987}.
Therefore~$\widehat g_r=(4-3r^2)g_r$ satisfies~$0<\sec_{\widehat g_r}\le1$.
Its volume is
\begin{align}
\operatorname{Vol}(\widehat g_r)
=(4-3r^2)^{(2m+1)/2}r\operatorname{Vol}(g_0)
\longrightarrow0.
\end{align}
Fix~$r$ so that this volume is smaller than the round volume.
From the scalar Weyl law~\cite{Chavel1984} we have
\begin{align}
\frac{\mu_j(\widehat g_r)}{\mu_j(g_0)}
\longrightarrow
\left(\frac{\operatorname{Vol}(g_0)}{\operatorname{Vol}(\widehat g_r)}\right)^{2/(2m+1)}>1.
\end{align}
By the signed Dirac Weyl law~\eqref{eq:round-signed-weyl}, the corresponding ratio for~$|\lambda_k|$ tends to the same volume ratio raised to the power~$1/(2m+1)$ as~$k\to\pm\infty$.
Thus both comparisons are strict at all sufficiently high indices of the fixed metric~$\widehat g_r$.
This argument does not imply that the first positive scalar eigenvalue diverges during collapse: functions pulled back from the quotient supply bounded eigenvalue branches.
\end{proof}

The construction uses the ordinary Hopf action in every odd dimension, since only the volume enters the Weyl asymptotics.

\subsection{The first Dirac gap and high-index comparison on even-dimensional spheres}
\label{subsec:even-spheres}

The following is a consequence of Klingenberg's injectivity-radius theorem and B\"ar's hyperspherical-radius estimate.
We reprove this deduction for completeness.

\begin{prop}\label{prop:even-upper}
Let~$(M^n,g)$ be a closed connected spin manifold of even dimension~$n\ge2$.
If~$0<\sec_g\le1$, then
\begin{align}\label{eq:even-gap-bound}
 \delta(g)\le\frac n2.
\end{align}
Equality holds if and only if~$(M,g)$ is the unit round sphere.
\end{prop}

\begin{proof}
The spin manifold is oriented and hence simply connected by Synge's theorem.
By Klingenberg's theorem~\cite{Klingenberg1959}, we have~$\operatorname{inj}(M,g)\ge\pi$.
Recall that the hyperspherical radius is defined by
\begin{align}\label{eq:hyperspherical-radius}
 \operatorname{Rad}_{\sph^n}(M,g)
 =\sup\left\{r>0:\begin{array}{l}
 \exists\,F:(M,g)\longrightarrow(\sph^n,r^2g_0),\\
 \operatorname{Lip}(F)\le1,\quad\deg F\ne0
 \end{array}\right\},
\end{align}
where~$g_0$ is the unit round metric.
By comparison on an exponential ball, we also have~$\operatorname{Rad}_{\sph^n}(M,g)\ge1$~\cite[Remark~7]{Baer2026Hyperspherical}.
B\"ar's theorem then yields
\begin{align}\label{eq:even-baer}
 \delta(g)^2\le
 \frac{n^2}{4\operatorname{Rad}_{\sph^n}(M,g)^2}\le\frac{n^2}{4}.
\end{align}
If~$\delta(g)=n/2$, both inequalities are equalities, and the rigidity statement in~\cite[Main Theorem]{Baer2026Hyperspherical} identifies~$M$ with the unit round sphere.
The converse follows from its Dirac spectrum.
Positive scalar curvature excludes harmonic spinors, so by chirality~$\delta(g)=\lambda_1(g)=|\lambda_{-1}(g)|$.
\end{proof}

B\"ar's estimate in terms of the hyperspherical radius holds in every dimension.
It is the lower bound for that radius under~$0<\sec_g\le1$ that uses evenness, through Synge's and Klingenberg's theorems.

\begin{rmk}\label{rmk:friedrich-round-class}
By Friedrich's estimate~\cite{Friedrich1980erste}, we have
\begin{align}\label{eq:friedrich-lower-comparison}
 |\lambda_k(g)|^2\ge\frac{n}{4(n-1)}\min_{\sph^n}R_g,
 \qquad\text{for each }k\in\bZ\setminus\{0\}.
\end{align}
Thus~$R_g\ge n(n-1)$ implies~$|\lambda_k(g)|\ge n/2$, and in particular the first eigenvalue of each sign is bounded away from zero by its round value.
This is a classical consequence and does not require a fixed conformal class.
If, in addition,~$g\in[g_0]$, where~$g_0$ is the unit round metric, equality for~$\lambda_1$ or~$\lambda_{-1}$ implies that~$g$ admits a real Killing spinor with Killing constant of absolute value~$1/2$.
The metric is then Einstein and conformally flat, hence isometric to the unit round sphere.
The same numerical lower bound does not compare~$\lambda_k(g)$ with the higher round levels~$|\lambda_k(g_0)|>n/2$.
For the upper direction, the use of Friedrich's twistor-spinor identities within~$[g_0]$ leads instead to the separate signed first-level bounds in Theorem~\ref{thm:round-conformal-first-level}.
\end{rmk}

\begin{cor}\label{cor:even-sphere-high-indices}
Let~$g_0$ be the unit round metric on~$\sph^n$, where~$n\ge2$ is even.
If~$0<\sec_g\le1$ and~$(\sph^n,g)$ is not isometric to the unit round sphere, then
\begin{align}
 |\lambda_k(g)|&<|\lambda_k(g_0)|,\quad\text{for all sufficiently large }|k|;\notag\\
 \mu_j(g)&<\mu_j(g_0),\quad\text{for all sufficiently large }j.\label{eq:lap-even-tail}
\end{align}
\end{cor}

\begin{proof}
Write~$v_g=\vol(\sph^n,g)$ and~$v_0=\vol(\sph^n,g_0)$.
By Klingenberg's injectivity-radius theorem and Berger's sharp isoembolic inequality~\cite{Klingenberg1959,Berger1980Isoembolic},
\begin{align}\label{eq:even-volume}
 v_g\ge v_0\left(\frac{\operatorname{inj}(\sph^n,g)}{\pi}\right)^n\ge v_0 .
\end{align}
Equality of the volumes forces~$\operatorname{inj}(\sph^n,g)=\pi$ and equality in Berger's inequality, so~$g$ is isometric to~$g_0$.
For a nonround metric we therefore have~$v_g>v_0$, and Remark~\ref{rmk:weyl-without-conformal-class} applies without any assumption on the conformal class of~$g$.
\end{proof}

\begin{rmk}\label{rmk:odd-even-parallel}
On spheres with~$0<\sec_g\le1$, the Weyl laws for the Dirac operator and the scalar Laplacian detect the same difference between odd and even dimensions.
In odd dimensions ordinary Hopf collapse produces metrics of volume smaller than the round volume, for which both spectra exceed their round counterparts at all sufficiently high indices.
In even dimensions Klingenberg's injectivity-radius bound and Berger's inequality reverse the volume comparison, and hence both spectral comparisons, unless the metric is round.
At the first Dirac eigenvalue it is the spin lift of the weighted circle action that matters in odd dimensions, whereas in even dimensions B\"ar's inequality converts the injectivity-radius bound into a bound for the first gap.
\end{rmk}

\

\textbf{Declaration of generative AI and AI-assisted technologies in the manuscript preparation process.}
The authors used ChatGPT (OpenAI) to obtain feedback on mathematical arguments and to improve the clarity and presentation of the manuscript.
The authors critically reviewed and revised the resulting material, independently verified the mathematical arguments and references, and accept full responsibility for the accuracy and integrity of the final manuscript.

\medskip

\textbf{Conflict of Interest} The authors have no conflicts to disclose. 

\medskip

\textbf{Data Availability Statement} Data sharing is not applicable to this article as no datasets were generated or analyzed during the current study.

\printbibliography[heading=bibintoc,title={References}]
\end{document}